\documentstyle[a4,11pt,amssymb,leqno]{amsart}

\title[Algebraic double reflection]
{Note on double reflection and algebraicity of holomorphic mappings}
\def\V{{\mathbb V}} 
\def\W{{\mathbb W}} 
\def\Y{{\mathbb Y}} 
\def\M{{\mathbb M}} 
\def\A{{\mathbb A}}
\def\L{{\mathbb L}} 
\def\X{{\mathbb X}}
\def\Y{{\mathbb Y}}
\def\N{{\mathbb N}} 
 
\def\S{{\mathbb S}} 
\def\Z{{\mathbb Z}} 
 
\def\R{{\mathbb R}} 
\def\C{{\mathbb C}} 
\def\L{{\mathbb L}}                                                          
\def\V{{\mathbb V}}                                                         
\date{\today}

\author{Jo\"el Merker}

\address{Laboratoire d'Analyse, Topologie et Probabilit\'es\\
Centre de Math\'ematiques et d'Informatique, UMR 6632\\
39 rue Joliot Curie\\
F-13453 Marseille Cedex 13, France}

\email{merker@@cmi.univ-mrs.fr, merker@@dmi.ens.fr}



\begin{document}

\begin{abstract} In this note, our purpose is to establish shortly
the algebraicity of a holomorphic mapping between real algebraic CR
manifolds under a double reflection condition which generalizes the
classical single reflection.  A complete study of various double
reflection conditions is also provided.
\end{abstract}

\maketitle

The goal of this note is to understand double reflection for
holomorphic mappings $f: M\to M'$ between real algebraic CR manifolds
in complex spaces of different dimensions. We plan to understand it in
the general case, that is with and without reducing, shrinking or
stratifying the first family of equations ``$f(z)\in r_{M'}
(f(Q_{\bar{z}})):= \V_z'$'' (see below) coming from first reflection.


Thus, let us quickly present the general problematic.

Let $f: M\to M'$ be a holomorphic map between real algebraic CR
manifolds in $\C^n$, $\C^{n'}$ given by real polynomial equations
$\rho(z, \bar{z})=0$, $\rho'(z',\bar{z}')=0$ so that $\rho'(f(z),
\bar{f}(\bar{w}))=0$ if $\rho(z,\bar{w})=0$ and let $Q_{\bar{w}}$,
$Q_{\bar{w}'}'$ denote Segre varieties.

It appears that two crucial observations yield heuristical insights
into the problem of finding conditions (C) such that: (C)
$\Rightarrow$ $f$ is algebraic.

\def\theequation{I}

First, in the observation that \begin{equation} f(z)\in
                      r_{M'}(f(Q_{\bar{z}})):=\{w' \ \!  {\bf :} \ \!
                      Q_{\bar{w}'}' \supset f(Q_{\bar{z}}) \}:= \V_z'
                      \end{equation} (we intentionally do not write
                      $w'\in U'$ or $w'\in \C^{n'}$: these two
                      possibilities will be studied hereafter) one is
                      lead to guess that the algebraic set $\V_z'$
                      which is parametrized by $z$ is a {\it finite
                      algebraic determinacy set} for the value of
                      $f(z)$ if ${\rm dim}_{f(z)} \V_z'=0$ $\forall \
                      z$. Being classical, this circumstance entails
                      that $f$ is algebraic, indeed ([SS], [BR],
                      [BER], [CPS]\footnote{ The set $\V_0'$ is called
                      the ``characteristic variety of $f$ at $0$'' in
                      [CPS].}).  As usual, this determinacy set
                      $\V_z'$ can be constructed by simply applying
                      the tangential Cauchy-Riemann operators of $M$
                      to the equations $\rho'(f(z),
                      \bar{f}(\bar{w}))=0$.

\def\theequation{II}

The second observation, which is due to Zaitsev [Z98], is
that\footnote{ We can call identity (I) ``first reflection'' and
identity (II) ``double reflection''.}  \begin{equation} f(z)\in
r_{M'}(f(Q_{\bar{z}})) \cap r_{M'}^2(f(Q_{\bar{w}})) := \V_z'\cap
r_{M'}(\V_w') := \X_{z,\bar{w}}' \end{equation} and one is lead to
analogously predict that $f$ is algebraic provided ${\rm dim}_{f(z)}
\X_{z,\bar{w}}'=0$ $\forall \ z, w$ such that
$\rho(z,\bar{w})=0$\footnote{ For reasons explained below, third
determination $r_{M'}^3$ or $r_{M'}^4$, {\it etc.}  provide no new
information.}.

(The choice of notation $\X_{z,\bar{w}}'$ instead of $\X_{z,w}'$ is
explained by the fact that the operator $r_{M'}$ conjugates the
dependence upon a parameter.)

\def\theequation{III}

Also, let us mention a possible third approach ([Z98], [D99]).  This
approach consists in choosing {\it smaller} algebraic sets $\W_z'
\subset \V_z'$ with suitable regularity properties\footnote{ Reducing
$\V_z'$ to being a {\it holomorphic family} with respect to $z$
[Z98].}  in order to compute the second reflection $r_{M'}(\V_w')$
easily. In particular, the shrinking from $\V_z'$ to $\W_z'$ is (and
in fact {\it must be}) constructive.  We will in this paper provide a
manner of constructing such a set $\W_z'$ in an uniform and
unambiguous way using minors of matrices of holomorphic functions
only.

This leads to a third type of determinacy: \begin{equation} f(z)\in
\W_z' \cap r_{M'} (\W_w') := \Z_{z,\bar{w}}'.  \end{equation}

Through various statements and examples, our aim will be thus to
compare conditions of determinacy of the value $f(z)$ by (I), (II) or
(III), to ask which one is weaker than another, to ask wether some are
necessary, to ask wether it is sufficient to require ${\rm dim}_{f(p)}
\V_p'=0$ or ${\rm dim}_{f(p)} \X_{p, \bar{p}}'=0$ or ${\rm dim}_{f(p)}
\Z_{p,\bar{p}}'=0$ for some or all $p\in M$, {\it etc.} Yet we will
mainly establish that:

1. Determination of $f(z)$ by $\W_z'$ is strictly finer than by
$\V_z'$\footnote{ We mean ${\rm dim}_{f(z)} \W_z' = 0$, ${\rm
dim}_{f(z)} \V_z' \geq 1$ $\forall \ z$ for some explicit examples of
$f$, $M$, $M'$, idem for 2, 3, 4, 5, 6 below.  These examples are
constructed in the paper.  }.

2. Determination of $f(z)$ by $\Z_{z,\bar{w}}'$ is strictly finer than
by $\W_z'$ (or $\V_z'$).

3. Determination of $f(z)$ by $\X_{z,\bar{w}}'$ is strictly finer than
   by $\V_z'$.

\noindent
{\it i.e.} that the inclusions $\X_{z,\bar{w}}' \subset \V_z'$ and
$\Z_{z,\bar{w}}'\subset \W_z'\subset \V_z'$ are all strict in general.

However, none of the inclusion $\Z_{z,\bar{w}}' \subset
\X_{z,\bar{w}}'$ and $\X_{z,\bar{w}}' \subset \Z_{z,\bar{w}}'$ is true
in general, because:

4. Determination of $f(z)$ by $\X_{z,\bar{w}}'$ can be strictly finer
than by $\Z_{z,\bar{w}}'$.

5. Determination of $f(z)$ by $\Z_{z,\bar{w}}'$ can be strictly finer
than by $\X_{z,\bar{w}}'$.

\def\theequation{IV}

\noindent
To sum-up, we are lead to define a fourth determinacy set
                      \begin{equation} f(z)\in \W_z' \cap
                      r_{M'}(\V_w') := \M_{z,\bar{w}}'=
                      \X_{z,\bar{w}}' \cap \Z_{z, \bar{w}}'.
                      \end{equation} and we have that:

6. Determination of $f(z)$ by $\M_{z,\bar{w}}'$ is strictly finer than
by $\X_{z,\bar{w}}'$ or $\Z_{z,\bar{w}}'$.

\noindent
{\it i.e.} that the inclusions $\M_{z,\bar{w}}' \subset
\X_{z,\bar{w}}'$ and $\M_{z,\bar{w}}' \subset \Z_{z,\bar{w}}'$ are all
strict in general.

Our goal is to find examples which exhibit the nonanalytic behavior of
$\X_{z,\bar{w}}'$ with respect to parameters and to explain why the
various second reflection conditions (II), (III) and (IV) are
inequivalent.

Our work deliberately forsakes the point of view of jet
parametrization about which the reader is referred to the works
[BER97], [Z97], [Z98], [D99], [ME99b].

This paper originates from questions the author have asked to Dmitri
Zaitsev at the Mathematische Forschungsinstitut (Oberwolfach,
Deutschland) during the Conference {\it Reelle Methoden der Komplexen
Analysis}, March 1998.

\section{Introduction and statement of results}

The general assumption throughout this article is:

$({\cal GH})$ $f: M\to M'$ is a local holomorphic map ${\cal
V}_{\C^n}(M)\footnote{ Notation: ${\cal V}_{\C^n}(M)$, ${\cal
V}_{\C^n}(p)$: small open neighborhood of $M$, of $p$ in $\C^n$, a
polydisc in case of a single point $p$.} \to \C^{n'}$, $M\to M'$,
between real algebraic CR generic manifolds in $\C^n$, $\C^{n'}$,
$p\in M$ is some point, $M$ is {\it minimal} at $p$ in the sense of
Tumanov, or equivalently of finite type in the sense of
Bloom-Graham. Even if part of our statements remain true for $M$
nonminimal, we shall not notify it for simplicity. Also, $n\geq 2$.

This general assumption is thought to be held throughout.

Let $m={\rm dim}_{CR} M$, $d={\rm codim}_{\R} M$, $m'={\rm dim}_{CR}
M'$, $d'={\rm codim}_{\R} M'$, $m+d=n$, $m'+d'=n'$, $m\geq 1$, $m'\geq
1$.

In suitable coordinates $z$, $z'$, $p=0$, $f(p)=0$, $M=\{z\in U \ \!
{\bf :} \ \!  \rho(z,\bar{z})=0\}$, $M'=\{z'\in U' \ \! {\bf :} \ \!
\rho'(z',\bar{z}')=0\}$, $U$, $U'$ are small polydiscs centered at the
origin, $\rho_j(z,\bar{z})=\sum_{|\mu|, |\nu| \leq N} \rho_{j,\mu,
\nu} z^{\mu} \! \bar{z}^{\nu}$, $1\leq j\leq d$,
$\rho_j'(z',\bar{z}')=\sum_{|\mu'|, |\nu'| \leq N'} \rho_{j,\mu',
\nu'} {z'}^{\mu'} \! \bar{z'}^{\nu'}$ $1\leq j\leq d'$ are {\it real
polynomials}, $\partial\rho_1 \wedge \cdots \wedge \partial
\rho_d(0)\neq 0$, $\partial\rho_1' \wedge \cdots \wedge \partial
\rho_{d'}'(0)\neq 0$.

We can assume $U=(\varepsilon \Delta)^n$, $\varepsilon >0$, $U' =
(\varepsilon ' \Delta)^{n'}$, $\varepsilon' > 0$, so $\overline{U}$
(conjugate) $=U$, $\overline{U}'= U'$.

Let $Q_{\bar{w}}=\{z\in U \ \! {\bf :} \ \!  \rho(z,\bar{w})=0\}$
denote Segre variety.

\def\theequation{1.\arabic{equation}} \setcounter{equation}{0}

For every subset $E\subset U$, define the action of the {\it
reflection operator}: \begin{equation} r_M(E):=\{w\in U\ \! {\bf :} \
\!  Q_{\bar{w}} \supset E\} =\footnotemark \cap_{w\in E} Q_{\bar{w}},
\ \ r_M^2(E) =r_M(r_M(E)), \end{equation} say, the {\it first
reflection} of $E$ and its {\it second reflection} across $M$.

Their basic properties are\footnotemark[7]:

1. $r_M(E)\cap E \subset M$.

2. $r_M(r_M(E)) \supset E$.

\noindent
Observe that $\cap_{k\in \N} r_M^k(E) = E\cap r_M(E)$. A similar
$r_{M'}$ is defined across $M'$.

Now, let $z\in Q_{\bar{w}}$. Let $f: M\to M'$ as in $({\cal GH})$.
Then $f(Q_{\bar{z}})\subset Q_{\overline{f(z)}}'$\footnotemark[7].
Also

1. $f(z) \in f(Q_{\bar{w}}) \subset r_{M'}^2(f(Q_{\bar{w}})).$

2. $f(Q_{\bar{z}}) \subset Q_{\overline{f(z)}}'$ hence by (1.1)
\begin{equation}
f(z)\in r_{M'}(f(Q_{\bar{z}})):= \V_z'.
\end{equation}

\footnotetext{See Lemma 2.1.}

 \noindent Consequently also, \begin{equation} f(z)\in
 r_{M'}(f(Q_{\bar{z}})) \cap r_{M'}^2(f(Q_{\bar{w}})) :=
 \X_{z,\bar{w}}'=\V_z'\cap r_{M'} (\V_{w}').  \end{equation}
 
 Observe that $r_{M'}^3$, $r_{M'}^4$, \ldots offer nothing more,
because $r_{M'}^{2k}(E')\supset E'$.

Also, a last notification:
\begin{equation}
r_{M'}(f(Q_{\bar{z}})) \subset Q_{\bar{f}(\bar{w})}', \ \ \ \ \
f(Q_{\bar{w}}) \subset r_{M'}^2(f(Q_{\bar{w}}))\subset
Q_{\bar{f}(\bar{w})}'.
\end{equation}

This introduction will now be divided in paragraphs a, b, c, d, e, f,
g, h, i, j, k, l, m corresponding to various questions and answers
that present themselves.

Organization of the paper consists in presentation of all the
problems, all the results, all the technical lemmas, all the examples
in the following paragraphs and of explanation of the all major links
between them. The checking of all the stated properties is reserved to
Sections 2, 3 and 4.

\smallskip
\noindent
{\bf a. The results.}

 The fundamental observation is that as $M'$ is algebraic, both
 $r_{M'}(E')$ and $r_{M'}^2(E')$ are {\it complex algebraic} sets for
 {\it any} set $E'$, since in (1.1) all the $Q_{\bar{w}'}'$ are.
 Therefore (1.3) should determine $f$ as an algebraic map of $z$ if
 ${\rm dim}_{f(z)} \X_{z,\bar{w}}'=0$ for all $z, w$ close to $0$,
 $z\in Q_{\bar{w}}$, a result which is true and which was originally
 proved in [Z98] for $\Z_{z,\bar{w}}'$ (and is well-known-classical
 with use of $\V_{z}':= r_{M'} (f(Q_{\bar{z}}))$ only\footnote{ The
 determination $f(z)\in r_{M'} (f(Q_{\bar{z}})) \cap
 Q_{\bar{f}(\bar{w})}'=\V_z'$ coincides with (1.2), since
 $r_{M'}(f(Q_{\bar{z}})) \subset Q_{\bar{f}(\bar{w})}'$.} instead of
 $\X_{z,\bar{w}}'$ or of $\Z_{z,\bar{w}}'$). See the closing remark
 p.25 here.

\smallskip
\noindent
{\bf Theorem 1.} {\it If ${\rm dim}_{f(z)} \V_z'=0$ or if ${\rm
dim}_{f(z)} \X_{z,\bar{w}}'=0$ $\forall \ z,w \in {\cal V}_{\C^n}(0)$,
$z\in Q_{\bar{w}}$, then $f$ is algebraic.}

\smallskip
{\it Remark.} Of course, the case ${\rm dim}_{f(z)} \V_z'=0$ is
contained in the more general case ${\rm dim}_{f(z)}
\X_{z,\bar{w}}'=0$, because clearly $\X_{z,\bar{w}}'\subset
\V_z'$. Our proof of Theorem 1 will be short using a partial
algebraicity theorem in [ME99a].

\smallskip
We shall indeed establish that the condition ${\rm dim}_{f(z)}
\X_{z,\bar{w}}'=0$, $\forall \ z, w\in {\cal V}_{\C^n}(0)$, $z\in
Q_{\bar{w}}$ implies that $f$ is complex algebraic on each Segre
variety in some neighborhood $V={\cal V}_{\C^n}(0)$ and then apply:

\smallskip
\noindent
{\bf Theorem 2.} ([ME99a]) {\it Let $g\in {\cal O}(V, \C)$, let $M$ be
minimal at $0$. Then $g$ is algebraic if and only if
$g|_{Q_{\bar{w}}\cap V}$ is algebraic $\forall \ w\in V$.}

\smallskip
We also obtain an equivalent version of Theorem 1:

\smallskip
\noindent
{\bf Theorem 1'.} ([Z98], [ME99a]) {\it If ${\rm dim}_{f(z)}
\X_{z,\bar{z}}'=0$ $\forall z\in M\cap {\cal V}_{\C^n}(0)$, then $f$
is algebraic.}

\smallskip
Theorem 1' admits several applications and covers several known
results ({\it e.g.}  [SS], [BR], [BER]), see [Z98]. Theorem 1' is the
original formulation of Theorem 1.1 in [Z98].  See the closing remark
p.24 here. In truth, some phenomena and subtle things are hidden
behind Theorem 1'.

Our work is aimed to reveal all of them.

\def\theequation{1.\arabic{equation}} \setcounter{equation}{3}

\smallskip
\noindent
{\bf b. First remarks and questions.}

This Theorem 1' will be deduced from Theorem 1 by proving that there
exist points $p\in M$ arbitrarily close to $0$ such that ${\rm
dim}_{f(z)}\X_{z,\bar{w}}' =0$ $\forall \ z, w \in {\cal
V}_{\C^n}(p)$, $z\in Q_{\bar{w}}$, see Proposition 5. One of the main
phenomenon here is that the set $\X_{z,\bar{w}}'$ is {\it not} in
general holomorphically parametrized by $z, \bar{w} \in U$, $z\in
Q_{\bar{w}}$, in the sense that there would exist analytic equations
such that $\X_{z,\bar{w}}'=\{z' \ \! {\bf :} \ \!
\lambda_j(z,\bar{w}, z')=0, 1\leq j\leq J\}$, such equations which
would readily imply that \begin{equation} {\rm dim}_{f(0)}\X_{0,0}' =0
\Rightarrow {\rm dim}_{f(z)} \X_{z,\bar{w}}' = 0, \ \forall \ z,w \in
{\cal V}_{\C^n}(0), z\in Q_{\bar{w}}.  \end{equation} And (1.4) above
would imply that Theorem 1 $\Rightarrow$ Theorem 1', but (1.4) fails.

Our goal is to explore properties of the map $(z,\bar{w})\mapsto
\X_{z,\bar{w}}'$. As a preliminary exposition, we will first recall
and state the well-behaved properties of the mapping $z\mapsto \V_z'$
right now in paragraph {\bf c} below.

Let us denote ${\cal M}=\{z\in Q_{\bar{w}}\}=\{(z,\bar{w}) \in U\times
U \ \! {\bf :} \ \!  \rho(z,\bar{w})=0\}$, which is a complex
$d$-codimensional submanifold of $U\times U$ called the {\it
complexification} of $M$.  We wish also to compare the twelve
conditions
$${\cal C}^4({\cal M}) \ : \ \ \ {\rm dim}_{f(z)}\M_{z,\bar{w}}'=0 \ \
\forall (z,\bar{w})\in {\cal M} \ \ \ \ \ C^4(M) \ : \ \ {\rm
dim}_{f(p)} \M_{p,\bar{p}}' = 0 \ \ \forall p\in M$$
$${\cal C}^3({\cal M}) \ : \ \ \ {\rm dim}_{f(z)}\X_{z,\bar{w}}'=0 \ \
\forall (z,\bar{w})\in {\cal M} \ \ \ \ \ C^3(M) \ : \ \ {\rm
dim}_{f(p)} \X_{p,\bar{p}}' = 0 \ \ \forall p\in M$$
$${\cal C}^2({\cal M}) \ : \ \ \ {\rm dim}_{f(z)}\Z_{z,\bar{w}}'=0 \ \
\forall (z,\bar{w})\in {\cal M} \ \ \ \ \ C^2(M) \ : \ \ {\rm
dim}_{f(p)} \Z_{p,\bar{p}}' = 0 \ \ \forall p\in M$$
$${\cal C}^1({\cal M}) \ : \ \ \ {\rm dim}_{f(z)}\V_{z,\bar{w}}'=0 \ \
\forall (z,\bar{w})\in {\cal M} \ \ \ \ \ C^1(M) \ : \ \ {\rm
dim}_{f(p)} \V_{p}' = 0 \ \ \forall p\in M$$ and also
$$C_p^4 \ : \ \ \ {\rm dim}_{f(p)} \M_{p,\bar{p}}' =0, \ \ \ \ C_p^3 \
: \ \ \ {\rm dim}_{f(p)} \X_{p,\bar{p}}' =0, \ \ \ \ \ \ \ \ \ \ \ \ \
\ \ \ $$
$$ \ \ \ \ \ \ \ \ \ \ \ \ \ \ \ \ C_p^2 \ : \ \ \ {\rm dim}_{f(p)}
\Z_{p,\bar{p}}' =0, \ \ \ \ C_p^1 \ : \ \ \ {\rm dim}_{f(p)} \V_{p}'
=0.$$

\smallskip
\noindent
{\bf c. Analytic dependence of $z\mapsto \V_z'$.}

Let ${\cal A}_n(U)= \{{\rm polynomials \ in} \ z\in U\}$, ${\cal
O}_n(U)=\{g\in {\cal O}(U, \C)\}$.  We abusively denote ${\cal
O}_n(U)\times {\cal A}_n(U)=\{g=g(z,w) \in {\cal O}_n(z)[w]\}$.  Also,
$\chi \in \overline{\cal O}_n(U)$ means $\chi$ antiholomorphic of a
variable $w\in U$: $\chi= \chi(\bar{w})$, $w\in U$.

By a well-known process of applying the tangential Cauchy-Riemann
operators to the identity $\rho'(f(z), \bar{f}(\bar{w}))=0$ as
$\rho(z,\bar{w})=0$, one can establish that the map $z\mapsto \V_z'$
is analytic (algebraic here because $M$, $M'$ are algebraic):

\smallskip
\noindent
{\bf Proposition 3.} {\it There exist $J\in \N_*$ and functions
$r_j(z,\bar{w}, z')\in {\cal A}_n(U)\times \overline{\cal O}_n(U)
\times {\cal A}_{n'}(U')$, $1\leq j\leq J$, such that $\forall \
(z,w)\in U\times U$, $z\in Q_{\bar{w}}$, \begin{equation}
r_{M'}(f(Q_{\bar{z}}))= \{z'\in U' \ \! {\bf :} \ \!  r_j(z,\bar{w},
z')=0, 1\leq j \leq J\}= \V_z':= \S_{z,\bar{w}}^1.  \end{equation}}

Here, we denote\footnote{By convention, we shall ``vectorially''
denote henceforth a set $r_1=0, \ldots, r_J:=0$, {\it i.e.}  simply by
$r=0$.}  \begin{equation} \S^1=\{(z,\bar{w}, z')\in U\times
\overline{U}\times U' \ \! {\bf :} \ \!  r(z,\bar{w}, z')=0\}
\end{equation} and $\S_{z,\bar{w}}^1$ denotes fibers of
$\S^1$. Although the equations of $\V_z'$ do depend in general of some
$\bar{w}$ such that $z\in Q_{\bar{w}}$, the zero-set $\V_z'$ appears
to be indeed independent of $\bar{w}$ provided $(z,\bar{w})\in {\cal
M}$ because it is by definition equal to $r_{M'}(f(Q_{\bar{z}}))$.
The equations $r(z,\bar{w}, z')$ justify the notation
$\S_{z,\bar{w}}^1$.

Proposition 3 and the upper semi-continuity of the fiber dimension of
a holomorphic map immediately imply that \begin{equation} C_0^1 \
\Rightarrow \ {\cal C}^1({\cal M}\cap (V\times V)) \ \ {\rm and} \ \
C_0^1 \Rightarrow C^1(M\cap V) \ \ {\rm for \ some} \ \ V={\cal
V}_{\C^n}(0)\subset U.  \end{equation}

\smallskip
\noindent
{\bf Corollary 4.} {\it If ${\rm dim}_{f(z)}\S_{z,\bar{w}}^1=0$ for
some $z\in U$ and some $w\in Q_{\bar{z}}$, then there exists ${\cal
N}\subset {\cal M}$ a proper complex analytic subvariety such that the
map \begin{equation} {\cal V}_{\C^{n'}}(f(z)) \ni z' \mapsto
r(z,\bar{w}, z') \in \C^J \end{equation} is an immersion at $f(z)$,
for $(z, \bar{w}) \in {\cal M} \backslash {\cal N}$.}

\smallskip
This is of course equivalent to the generic rank of the mapping
\begin{equation} {\cal M} \times \C^{n'} \ni (z,\bar{w}, z') \mapsto
(z,\bar{w}, r(z,\bar{w}, z')) \in {\cal M} \times \C^J \end{equation}
being equal to $2m+d+n'$. Just one remark more. As $M\cong
\{(z,\bar{z}) \in {\cal M}\}$ is a maximally real, real algebraic
submanifold of ${\cal M}$, then ${\cal N} \cap M:= N$ is a proper real
analytic subset of $M$ also.  In particular, after applying the
implicit function theorem to (1.8) near $(p, \bar{p})$, we obtain the
existence of $V_p={\cal V}_{\C^n}(p)$ and of
$\Psi_{\nu}'(z,\bar{w})\in {\cal A}_n (V_p)\times \overline{\cal
O}_n(V_p)$, $1\leq \nu \leq n'$ such that $f_{\nu}(z) =
\Psi_{\nu}'(z,\bar{w})$, $\nu= 1,\ldots,n'$, $z\in Q_{\bar{w}}$ (that
$\Psi_{\nu}'$ can be supposed polynomial in $z$ is achieved through
elimination theory).  Fixing $w$, this shows that $f$ is algebraic on
Segre varieties and Theorem 2 applies to show that $f$ is algebraic.

For further properties and knowledge about the geometry of $\S^1$ (the
first reflection variety), we refer to [BR], [BER], [BER97], [CMS],
[Z98], [BER99], [CPS], [ME99a].

\smallskip
\noindent
{\bf d. Almost everywhere analytic dependence of $(z,\bar{w})\mapsto
\X_{z,\bar{w}}'$.}

Now, we present the way how $(z,\bar{w})\mapsto \X_{z,\bar{w}}'$
varies:

\smallskip
\noindent
{\bf Proposition 5.} {\it If ${\rm dim}_{f(z)}\X_{z,\bar{z}}'=0$,
$\forall \ z\in M$, then there exists a dense open subset $D_M\subset
M$ such that

$(*)$ $\forall \ p\in D_M$, $\exists U_p={\cal V}_{\C^n}(p)$, $\exists
U_{p'}' = {\cal V}_{\C^{n'}}(p'=f(p))$, $\exists K\in \N_*$, $\exists
(s_k)_{1\leq k \leq K}$, $s_k\in \overline{\cal O}_n({U}_p) \times
{\cal O}_n(U_p) \times {\cal A}_{n'}(U_{p'}')$ such that, if
$\S^1=\{r(z, \bar{w}, z')=0\}$,

$\forall \ z,w \in U_p$, $z\in Q_{\bar{w}}$, $\forall \ w, z_1\in
U_p$, $w\in Q_{\bar{z}_1}$, \begin{equation} U_{p'}'\cap
\X_{z,\bar{w}}'\subset \{z'\in U_{p'}' \ \! {\bf :} \ \!  r(z,\bar{w},
z') =0, s(\bar{w}, z_1, z')=0\}:= \S_{z,\bar{w}, z_1}^2 \end{equation}
and such that the graph of $f$ over ${\cal M}\sharp {\cal M}:=\{
(z,\bar{w}, z_1) : \rho(z,\bar{w})=0, \rho(w,\bar{z}_1)=0\}$
intersected with $U_p\times \overline{U}_p \times U_p$ satisfies
\begin{equation} \Gamma r (f) =\{(z, \bar{w}, z_1, f(z))\}= \{(z,
\bar{w}, z_1, z') \ \! {\bf :} \ \!  z'\in \S_{z,\bar{w}, z_1}^2\}:=
\S^2.  \end{equation} Furthermore, there exist similar analytic
equations $r(z, \bar{w}, z') = 0$, $s(\bar{w}, z_1, z') = 0$ such that

$(**)$ $\forall \ (z,\bar{w}, z_1) \in ({\cal M} \sharp {\cal M}) \cap
(U_p\times \overline{U}_p \times U_p)$, the map \begin{equation} {\cal
V}_{\C^{n'}}(f(z)) \ni z' \mapsto (r(z,\bar{w}, z'), s(\bar{w}, z_1,
z'))\in \C^{J+K} \end{equation} is an immersion at $f(z)$.}

\smallskip
We invite the reader to notice that dependence of $\S^2$ is
holomorphic in $z$, antiholomorphic in $w$, which justifies and
explains the notation $\X_{z,\bar{w}}'$.

This proposition appeals several remarks.

The first one is: what is the structure of the closed set $M\backslash
D_M \subset M$ ?  Leaving this question for a while, we return to it
in Examples 13 and 14.

The second one is: we prove Theorem 1 without using Proposition 5.
This proposition is indeed used only to prove that Theorem 1
$\Rightarrow$ Theorem 1'.

Next, a third remark.  As $\Gamma r (f) \equiv \S^2$, the projection
$\pi : \S^2 (\subset U_p \times \overline{U}_p \times U_p \times
U_{p'}') \to U_p \times \overline{U}_p \times U_p$ is submersive.  The
zero locus $\S^2 \equiv \Gamma r (f) $ is smooth, but the equations
defining $\S^2$ can be nonreduced. After taking the reduced complex
space ${\rm Red}\ \S^2$, we obtain $(**)$. Finally, we obtain the
immersion property in $(**)$ of Proposition 5 for $z, w, z_1 \in {\cal
V}_{\C^n} (p)$.

Now, the main remark. Since clearly \begin{equation} {\rm
                      dim}_{f(p)}\S_{p,\bar{p}, p}^2 = 0 \ \Rightarrow
                      \ {\rm dim}_{f(z)}\S_{z,\bar{w}, z_1}^2 = 0 \
                      \forall z,w, z_1 \in W_p = {\cal V}_{\C^n}(p)
                      \subset U_p, \end{equation} we obtain as desired
                      that \begin{equation} {\rm
                      dim}_{f(z)}\X_{z,\bar{z}}'=0 \ \forall \ z\in M
                      \ \Rightarrow \ {\rm dim}_{f(z)}\X_{z,\bar{w}}'
                      = 0 \ \forall \ z,w \in W_p, z\in Q_{\bar{w}},
                      \end{equation} whence the reduction of Theorem
                      1' to Theorem 1 through Proposition 5 is
                      completed. This technical proposition is
                      interpreted in geometric terms of ``holomorphic
                      families'' in [Z98].

Finally, the equations $s(\bar{w}, z_1, z')$ of which Proposition 5
asserts the existence are clearly those for
$r_{M'}^2(f(Q_{\bar{w}}))$, while as before $r(z, \bar{w}, z')$ come
for $r_{M'}(f(Q_{\bar{z}}))$.

\smallskip
\noindent
{\bf e. Solvability of $f$ over a dense open set.}

The fundamental remark is that after solving from $(**)$ of
Proposition 5 at $q\in U_p$, $p\in D_M$ the collection of equations
\begin{equation} r_j(z,\bar{w}, f(z))=0, \ \ 1\leq j\leq J, \ \
s_k(\bar{w}, z_1, f(z))=0, \ \ 1\leq k\leq K, \end{equation} where
$z\in Q_{\bar{w}}$, $w\in Q_{\bar{z}_1}$, we have:

\smallskip
\noindent
{\bf Corollary 6.} {\it $\forall \ p\in D_M$, $\forall \ q\in U_p=
{\cal V}_{\C^n}(p)$, $\exists W_q= {\cal V}_{\C^n}(q)$, $\exists
\Psi_{\nu}'(z,\bar{w}, z_1)\in {\cal A}_n(W_q)\times \overline{\cal
O}_n({W}_q) \times {\cal O}_n(W_q)$, $1\leq \nu \leq n'$, such that
\begin{equation} f_1(z)=\Psi_1'(z,\bar{w}, z_1), \ldots ,
f_{n'}(z)=\Psi_{n'}'(z,\bar{w}, z_1), \end{equation} $\forall z,
w,z_1\in W_q$, $z\in Q_{\bar{w}}$, $w\in Q_{\bar{z}_1}$.}

\smallskip
\noindent
and equation (1.16) immediately explains that $f$ is algebraic on each
Segre variety $Q_{\bar{w}} \cap W_q$: just fix $\bar{w}, z_1$ in
(1.16) and let $z\in Q_{\bar{w}}$ vary ({\it cf.} [D99]).

\smallskip
\noindent
{\bf f. Algebraicity of $f$.}

Finally, Theorem 1' admits the following main corollary:

\smallskip
\noindent
{ \bf Theorem 7.} ([CMS], [Z98], [ME99a]) {\it If $M'$ does not
contain complex algebraic sets of positive dimension, then $f$ is
algebraic.}

\smallskip
{\it Proof.}  Indeed,

1. $\X_{z,\bar{z}}'=r_{M'}(f(Q_{\bar{z}}))\cap
r_{M'}^2(f(Q_{\bar{z}}))\subset M'$ (by $r_{M'}(E') \cap r_{M'}^2(E')
\subset M'$ $\forall \ E'$).

2. $r_{M'}(f(Q_{\bar{z}}))$, $ r_{M'}^2(f(Q_{\bar{z}}))$,
$\X_{z,\bar{z}}'$ are complex algebraic sets through $f(z)$.

\noindent
so that ${\rm dim}_{f(z)} \X_{z,\bar{z}}' =0$ $\forall \ z\in M$
necessarily in the assumption of Theorem 3: Theorem 1' applies. \hfill
$\square$

\smallskip
{\it Remark.} The author obtains a completely different proof of
Theorem 7 in [ME99a]. The proof in [CMS] is obtained for $M$
Segre-transversal instead of being minimal.

\smallskip
We have now completed the presentation of the main steps in the proof
of Theorem 1'.

\smallskip
\noindent
{\bf g. Comparison of $\V_z'$, $\X_{z,\bar{w}}'$.}

Next, we come to the comparison between conditions about $\V_z'$ and
$\X_{z,\bar{w}}'$.  It is known that:

There exist $f$, $M$, $M'$, $U$, $U'$ such that

1. ${\rm dim}_{f(z)}\V_z' \geq 1$, $\forall \ z\in U$ and

2. ${\rm dim}_{f(z)}\X_{z,\bar{w}}'=0$, $\forall \ z,w\in U$, $z\in
   Q_{\bar{w}}$,

\noindent
whence ${\cal C}^3({\cal M})$ is strictly finer that ${\cal C}^1({\cal
M})$ (identical fact about ${\cal C}^2({\cal M})$ or ${\cal C}^4({\cal
M})$).

\smallskip
\noindent
{\bf Example 8.} Take\footnote{Details are left to the reader. See
also further examples below.} $M: z_4=\bar{z}_4 + iz_1\bar{z}_1$ in
$\C^2_{(z_1, z_4)}$, $f(z_1, z_4)=(z_1, 0, 0, z_4)\in \C^4$, and the
hypersurface \begin{equation} M': \ \ z_4'=\bar{z}_4' +
iz_1'\bar{z}_1'+i{z_1'}^2 \bar{z}_3'+
i{\bar{z'}_1}^2z_3'+iz_2'\bar{z}_3'+ i\bar{z}_2' z_3'.  \end{equation}

{\it Remark.} Second reflection is superfluous in case $n=n'$, $f$ is
a biholomorphic map $U\to U'$ ({\it cf.} [BR], [BER], [BR97], [Z98])
or if one assumes directly that ${\rm dim}_{f(z)} \V_z'=0$ ({\it cf.}
[CPS]).

\smallskip
\noindent
{\bf h. Comparison of $\V_z', \W_z', \Z_{z,\bar{w}}',
\X_{z,\bar{w}}'$.}

 Yet another strategy ({\it cf.} [Z98], [D99]) consists in replacing
if possible the set $\S^1=\{(z,\bar{w}, z') \ \! {\bf :} \ \!
r(z,\bar{w}, z') =0\}$ by some smaller complex analytic set
$\widetilde{\S}^1=\{(z,\bar{w}, z')\ \! {\bf :} \ \!
\tilde{r}(z,\bar{w}, z')=0\}\subset \S^1$ such that

1. $\tilde{r} \in {\cal O}_n(U) \times \overline{\cal O}_n({U}) \times
{\cal A}_{n'}(U')$.

2. $\Gamma r (f)=$ the graph of $f$ $=\{(z,\bar{w}, f(z)) \ \! {\bf :}
\ \!  (z,\bar{w})\in {\cal M}\} \subset \widetilde{\S}^1$.

3. $\widetilde{\S}^1$ is obtained in a constructive way.

\noindent
That $\widetilde{\S}^1$ should be given by means of an explicit
construction is important, because the datum is $\S^1$ from which one
tries to deduce that $z'$ is solvable in terms of $z, \bar{w}$.

Taking such a set $\widetilde{\S}^1$, one can expect that ${\rm
dim}_{f(z)}\widetilde{\S}^1_{z,\bar{w}}=0$. For instance, if $\Gamma
r(f)$ is contained in ${\rm Sing} (\S^1)$ which is computable in terms
of $r(z,\bar{w}, z')$ only since it is explicitely given, and because
$\Gamma r (f) \subset {\rm Sing} (\S^1)$ can be tested, it is possible
to shrink $\S^1$ and to replace it by $\widetilde{\S}^1:= {\rm Sing}
(\S^1)$, obtaining {\it new}, {\it possibly finer} equations
$\tilde{r}(z,\bar{w}, f(z))=0$ ({\it cf.} Lemma 4.2 in [Z98]). In
[Z98], $\S^1$ is also shrunk more again, still in a constructive way,
in order that $\widetilde{\S}^1$ becomes a {\it holomorphic
family}. Therefore, there might exist many different
$\widetilde{\S}^1$ depending on the way how $\S^1$ is shrunk in a
constructive way.  Uniqueness of $\widetilde{\S}^1$ is not clear in
[Z98].

However in the end of Section 3 we propose a uniform unambiguous
method which even uses only elementary tools: minors and the
uniqueness principle, and not passing to the filtration by singular
subspaces.

If ${\rm dim}_{f(z)}\widetilde{\S}^1_{z,\bar{w}} \geq 1$, denote
$\W_{z,\bar{w}}':= \{z'\in U'\ \! {\bf :} \ \!  \tilde{r}(z,\bar{w},
z')=0\}=\widetilde{\S}_{z,\bar{w}}^1$ (to recover the notations of our
Presentation temporarily) and \begin{equation}
\Z_{z,\bar{w}}':=\W_{z,\bar{w}}' \cap r_{M'} (\W_{w,\bar{z}_1}'), \ \
\ \ \ z\in Q_{\bar{w}}, \ w\in Q_{\bar{z}_1}.  \end{equation} Then
$f(z)\in \Z_{z,\bar{w}}'$ because

\noindent
(1.19) \hspace{1.9cm} $f(z)\in \W_{z,\bar{w}}' \subset \V_z'$

\noindent
\hspace{3cm} $r_{M'}(\W_{z,\bar{w}}') \supset r_{M'}(\V_z')$ (by
$r_{M'}(E') \supset r_{M'}(F')$ if $E' \subset F'$)

\noindent
\hspace{3cm} $r_{M'}(\W_{w,\bar{z}_1}') \supset
r_{M'}(\V_w')=r_{M'}^2(f(Q_{\bar{w}})) \supset f(Q_{\bar{w}}) \ni
f(z)$

\noindent
\hspace{3cm} $f(z)\in \W_{z,\bar{w}}' \cap r_{M'}(\W_{w,\bar{z}_1}') =
\Z_{z,\bar{w}, z_1}'$.

\def\theequation{1.\arabic{equation}} \setcounter{equation}{19}

\noindent
The gain in reducing $\S^1 \mapsto \widetilde{\S}^1$ lies in the fact
that one can easily insure that $(w,\bar{z}_1) \mapsto
r_{M'}(\W_{w,\bar{z}_1}')$ becomes an analytic parametrization
(holomorphic family, [Z98]) by having first a nice representation of
$\W_{z,\bar{w}}'$:

\smallskip
\noindent
{\bf Proposition 9.} ([Z98], [D99]) {\it There exists ${\cal N}
\subset {\cal M}$ a proper complex analytic subset such that $\forall
\ p=(z_p, \bar{w}_p) \in {\cal M} \backslash {\cal N}$ $\exists \
{\cal U}_p = {\cal V}_{\cal M}(p)$, $\exists \ n_1'$, $n_2'\in \N$,
$n_1'+n_2'=n'$, $\exists \ \Phi_{\nu}'(z,\bar{w}, z') \in {\cal
A}_n({z})\times {\cal O}_n(\bar{w}) \times {\cal A}_{n'}(z_1')$, $(z,
\bar{w})\in {\cal U}_p$, $(w, \bar{z}_1) \in {\cal U}_p$, $1\leq \nu
\leq n'_1$, such that \begin{equation} \{(z,\bar{w}, f(z)) \ \! {\bf
:} \ \!  (z,\bar{w})\in {\cal M} \} \subset \{(z,\bar{w}, z') \ \!
{\bf :} \ \!  z_2' =\Phi'(z,\bar{w}, z_1')\} \subset \{r(z,\bar{w},
z')=0\}.  \end{equation}}

In [D99], it is established that representation (1.20) is unique: the
set $\{z_2'=\Phi'(z, \bar{w}, z_1')\}$ being the maximal for inclusion
among the sets in the form $\Lambda=\{ z_2'=\Psi'(z, \bar{w}, z_1')\}$
(for some splitting of the coordinates $z'$) satisfying $\Gamma r(f)
\subset \Lambda \subset \S^1$.

Let $p\in M\backslash (N:={\cal N}\cap M)$, let $U_p={\cal
V}_{\C^n}(p)$ such that ${\cal U}_p \subset U_p \times U_p$.  This
representation yields after easy work that there exist $K\in \N$, $ \
(s_k)_{1\leq k\leq K} \in \overline{\cal O}_n ({U}_p) \times {\cal
O}_n(U_p) \times {\cal A}_{n'} (U_{p'}')$ such that $\forall \ z, w\in
U_p$, $z\in Q_{\bar{w}}$, $\forall \ w, z_1 \in U_p$, $w\in
Q_{\bar{z}_1}$, \begin{equation} \Z_{z,\bar{w}, z_1}'=\{z'\in U_{p'}'
\ \! {\bf :} \ \!  r(z,\bar{w}, z')=0, \ s(\bar{w}, z_1, z')=0\}.
\end{equation} Whence the holomorphicity of the map
$(z,\bar{w})\mapsto \Z_{z,\bar{w}}'$ and the set $\Z_{z,\bar{w},
z_1}'$ does not depend as a set of $z_1$ if $(z,\bar{w}, z_1)\in {\cal
M} \sharp {\cal M}$ and coincides with $\Z_{z,\bar{w}}'$ which was
defined in a set theoretical way.  See also Proposition 16.

\smallskip
\noindent
{\bf Fundamental remark.} The constructiveness of a shrinking
$\S^1\mapsto \widetilde{\S}^1$ is essential. One is temptated to
introduce $\S_{min}^1:=$ the minimal (for inclusion) ${\cal A}_n\times
\overline{\cal O}_n \times {\cal A}_{n'}$-set contained in $\S^1$
satisfying $\Gamma r (f) \subset \S_{min}^1\subset \S^1$, {\it i.e.}
the intersection of all $\widetilde{\S}^1$, even those not
constructive, and to put ${\Z_{min,z,\bar{w},z_1}}':= \S_{min,
z,\bar{w}}^1 \cap r_{M'}(\S_{min, w,\bar{z}_1}^1)$. However, the
equations $r_{min}(z,\bar{w}, z')=0$ being not known from the datum
$\S^1$ in general and not constructible in some explicit way, it is
quite impossible to deduce from $f(z)\in {\Z_{min,z,\bar{w}, z_1}}'$
anything. Not to mention that anyway if $f$ was algebraic from the
beginning, $r_{min}:= z'-f(z)$ would have been convenient and the
condition ${{\rm dim}_{f(z)}\Z_{min,z,\bar{w}, z_1}}'=0$ (here ${\rm
dim}_{f(z)} \S_{min, z,\bar{w}, z_1}^1=0$) then becomes tautological!

\smallskip 
Before entering in further discussions, let us summarize the
properties of $\Z_{z,\bar{w}}'$ as follows.  We recommand to see also
Proposition 16.

\smallskip
\noindent
{\bf Theorem 10.} {\it $(1)$ The set of points where
$(z,\bar{w})\mapsto \Z_{z, \bar{w}}'$ is not holomorphic is a proper
complex analytic subset ${\cal N}$ of ${\cal M}$.  Let $N:= {\cal N}
\cap M$.

$(2)$ If ${\rm dim}_{f(z)}\Z_{z,\bar{z}}' =0$ $\forall \ z\in {\cal
V}_{\C^n}(0)\cap M$, then ${\rm dim}_{f(z)}\Z_{z,\bar{z}}'=0$ for $z$
outside a proper real analytic subvariety $N_1$ of $M\backslash N$.

$(3)$ If ${\rm dim}_{f(p)} \Z_{p, \bar{p}}'=0$ at some point $p\in
M\backslash (N\cup N_1)$, then $f$ is algebraic.}

\smallskip
{\it Remark.} That the bad set ${\cal N}$ is analytic is a property
which is specific to $\Z_{z,\bar{w}}'$. For $\X_{z,\bar{w}}'$, the set
$M\backslash D_M$ is definitely not analytic, see Example 14.

\smallskip

We can now summarize the result and the actual proof given in [Z98].
During the course of the proof of Theorem 1' (Theorem 1.1 there), the
author shows that by shrinking $\S^1$ to $\widetilde{\S}^1$, on can
obtain the nice representation (1.20) by analytic equations. The place
is in Lemma 4.2 [Z98], where $\S^1$ is stratified. One question which
is left unsolved is to ask wether the study of $\Z_{z,\bar{w}}'$ and
of $\X_{z,\bar{w}}'$ are equivalent. What is actually proved is not
Theorem 1' but:

\smallskip
\noindent
{\bf Theorem.} ([Z98]) {\it If ${\rm dim}_{f(p)} \Z_{p,\bar{p}}'=0$
$\forall \ p\in {\cal V}_{\C^n}(0) \cap M$, then $f$ is algebraic.}

\smallskip
We nevertheless remark that this theorem implies anyway Theorem 7
(Theorem 1.3 in [Z98]) and that the proof of Theorem 7 from the above
theorem goes on identically, because $\Z_{p,\bar{p}}'= \W_{p,\bar{p}}'
\cap r_{M'} (\W_{p,\bar{p}}')\subset M'$ and is algebraic.

See the end of Section 3 for construction of $\widetilde{\S}^1$.

\smallskip
Now, we return to comparison of $\X_{z,\bar{w}}'$ with
$\Z_{z,\bar{w}}'$. If ${\rm dim}_{f(z)} \W_{z,\bar{w}}' \geq 1$
$\forall \ z \in U$ (so second reflection is needed), one can expect
that \begin{equation} {\rm dim}_{f(z)}\X_{z,\bar{w}}' = 0 \ \forall \
z \ \Rightarrow \ {\rm dim}_{f(z)} \Z_{z,\bar{w}}' = 0 \ \forall \ z,
\end{equation} {\it i.e.} that the study of $\Z_{z,\bar{w}}'$ is
sufficient to get a complete proof of Theorem 1.  Nevertheless, {\bf
two inclusions in the opposite side enter in competition}
\begin{equation} \W_{z,\bar{w}}' \subset \S_{z,\bar{w}}^1 \ \ \ \ \
{\rm and} \ \ \ \ \ r_{M'}(\W_{z,\bar{w}}') \supset
r_{M'}(\S_{z,\bar{w}}^1) \end{equation} so that it is not clear how
\begin{equation} \X_{z,\bar{w}, z_1}'=\S_{z,\bar{w}}^1 \cap
r_{M'}(\S_{w,\bar{z}_1}^1) \ \ \ \ \ {\rm and} \ \ \ \ \
\Z_{z,\bar{w}, z_1}'= \W_{z,\bar{w}}' \cap r_{M'} (\W_{w, \bar{z}_1}')
\end{equation} could be comparable. Indeed, implication (1.21) is
untrue:

\smallskip
\noindent
{\bf Example 11.} There exist $f$, $M$, $M'$ such that
                      \begin{equation} {\rm dim}_{f(z)} \X_{z,\bar{w},
                      z_1}' = 0 \ \ \forall \ z \ \ \ \ \ \hbox{but} \
                      \ \ \ \ {\rm dim}_{f(z)}\Z_{z,\bar{w}, z_1}'
                      \geq 1 \ \ \forall \ z.
\end{equation}
Explicitely, take: $M: z_5=\bar{z}_5+iz_1\bar{z}_ 1$ in $\C^2_{(z_1,
z_5)}$, $f(z_1, z_5)=(z_1,0,0,0, z_5)$, and \begin{equation} M': \ \
z_5'=\bar{z}_5' + iz_1'\bar{z}_1' +i {z_1'}^2 \bar{z}_3' \bar{z}_4'+
i{\bar{z'}_1}^2z_3'z_4'+i{z_3'}^2
{\bar{z'}_2}^2+i{\bar{z'}_3}^2{z_2'}^2+
\end{equation}
$$ \ \ \ \ \ \ \ \ \ \ \ \ \ \ \ \ \ \ \ \ \ \ \ \ \ \ \ \ \ \
+i{\bar{z'}_3}^3 {z'}_4 ^3+i{z'}_3^3 \bar{z'}_4^3.
$$
In conclusion, determination of $f(z)$ by $\X_{z,\bar{w}, z_1}'$ can
be strictly finer than by $\Z_{z,\bar{w}, z_1}'$.

For the details, see Section 4.

\smallskip
And very surprisingly, it is also true that determination of $f(z)$ by
$\Z_{z,\bar{w}}'$ can be strictly finer than by $\X_{z,\bar{w}}'$.

\smallskip
\noindent
{\bf Example 12.} There exist $f$, $M$, $M'$ such that ${\rm
dim}_{f(z)}\Z_{z,\bar{w}}'= 0$ $\forall \ z$ but ${\rm
dim}_{f(z)}\X_{z,\bar{w}}' \geq 1$ $\forall \ z$. To be explicit, take
$M: z_4=\bar{z}_4+ iz_1\bar{z}_1$ in $\C_{(z_1, z_4)}$, $f(z_1, z_4)=
(z_1,0,0,z_4) \in \C^4$, and \begin{equation} M': \ \ \ \ \
z_4'=\bar{z}_4'+iz_1'\bar{z}_1'+ i{z_1'}^2{\bar{z'}_2}\bar{z}_3 +
i{\bar{z'}_1}^2 z_2' z_3'.  \end{equation} In summary \begin{equation}
{\rm dim}_{f(z)}\X_{z,\bar{w}}'= 0 \ \forall \ z \ \not\Leftarrow \
\not\Rightarrow \ {\rm dim}_{f(z)}\Z_{z,\bar{w}}' = 0 \ \forall \ z
\end{equation} and therefore, it is justified to introduce
\begin{equation} \M_{z,\bar{w}}':=\W_{z,\bar{w}}' \cap r_{M'}(\V_w')
\end{equation} where \begin{equation} \W_{z,\bar{w}}' =
\widetilde{\S}^1_{z,\bar{w}} \end{equation} for a constructive
shrinking $\widetilde{\S}^1$ of $\S^1$.  (Of course, different such
shrinkings may exist, depending on the conditions that are
imposed. The choice $\widetilde{\S}^1=\S^1$ can always be done. Our
examples illustrate well the phenomenon.)

\smallskip
\noindent
{\bf i. Comparison of $\M_{z,\bar{w}}'$ and $\Z_{z,\bar{w}}'$,
$\X_{z,\bar{w}}'$.}

Notice that $\M_{z,\bar{w}}' = \X_{z,\bar{w}}'\cap \Z_{z,\bar{w}}'$.

Now it is clear that \begin{equation} {\rm dim}_{f(z)}\M_{z,\bar{w}}'
                      = 0 \ \forall \ z,w , z\in Q_{\bar{w}} \
                      \Leftarrow \ {\rm dim}_{f(z)}\X_{z,\bar{w}} ' =
                      0 \ \forall z, w, z\in Q_{\bar{w}}
                      \end{equation} \begin{equation} {\rm
                      dim}_{f(z)}\M_{z,\bar{w}}' = 0 \ \forall \ z,w ,
                      z\in Q_{\bar{w}} \ \Leftarrow \ {\rm
                      dim}_{f(z)}\Z_{z,\bar{w}} ' = 0 \ \forall z, w,
                      z\in Q_{\bar{w}} \end{equation} Our examples
                      also show that the reverse implications
                      $\Rightarrow$ are both untrue.

\smallskip
\noindent
{\bf j. Summarizing tabulae.}

It is time to give a complete link tabular between the twelve
conditions
                  $$
C_p^1, C_p^2, C_p^3, C_p^4, C^1(M), C^2(M), C^3(M), C^4(M), {\cal
C}^1({\cal M}), {\cal C}^2({\cal M}), {\cal C}^3({\cal M}), {\cal
C}^4({\cal M}).$$ Return to definitions. Here, $p\in M$ is a fixed
chosen point, the origin in previous coordinates. The point $p$ is
chosen arbitrarily and is fixed. $C_p^j$ denotes: $C_p^4: {\rm
dim}_{f(p)} \M_{p,\bar{p}}'=0$, $C_p^3: {\rm dim}_{f(p)}
\X_{p,\bar{p}}'=0$, $C_p^2: {\rm dim}_{f(p)} \Z_{p,\bar{p}}'=0$,
$C_p^1: {\rm dim}_{f(p)} \V_{p}'=0$. Here, $C^j(M)$ does not denote
``$C_q^j$ $\forall \ q \in M$, but $\forall \ q\in D_M$'' over a dense
open subset $D_M$ of $M$. Idem for ${\cal C}^j({\cal M})$. Therefore a
priori, the implication $C^j(M)\Rightarrow C_p^j$ is false, since $p$
can belong to the set of points $q$ where $C_q^j$ is not satisfied.

First, we know that ${\cal C}^j({\cal M}) \Longleftrightarrow_{\cal D}
C^j(M)$ modulo open dense sets, {\it i.e.} if $C^j(M)$ is satisfied
over an open dense subset of $M$, then ${\cal C}^j({\cal M})$ is
satisfied over an open dense subset of ${\cal M}$, $j=1,2,3,4$, and
conversely.  We know this because of Propositions 3, 5, 9.  Indeed,
for $j=1, 3$, see Propositions 3, 5. For $j=2, 4$, see Propositions 5,
9 and eq. (1.21) and also Proposition 1.16 below.  Therefore, if $C_p$
denotes $(C_p^j)_{1\leq j\leq 4}$, $C(M)=(C^j(M))_{1 \leq j\leq 4}$,
${\cal C}({\cal M})= ({\cal C}^j({\cal M}))_{1\leq j\leq 4}$, the
comparison of our twelve conditions which could have been explained in
a $12\times 12$ tabular with 144 entries can be reduced to three
$4\times 4$ tabulars:

{\scriptsize
$$
\begin{tabular}{||c||c||c||c||}
\hline\hline $***$ & $C_p$ & $C(M)$ & ${\cal C}({\cal M})$ \\
\hline\hline $C_p$ & $\times$ & $\times$ & $\times$ \\ \hline\hline
$C(M)$ & $\times$ & $\times$ & $\times$ \\ \hline\hline ${\cal
C}({\cal M})$ & $\times$ & $\times$ & $\times$ \\ \hline \hline
\end{tabular} \ \ \
\cong \ \ \ \left( \ \
\begin{tabular}{||c||c||}
\hline\hline $*$ & $C_p$\\ \hline\hline $C_p$ & $\times$ \\
\hline\hline
\end{tabular}, \ \ \ 
\begin{tabular}{||c||c||}
\hline\hline $*$ & $C(M)$\\ \hline\hline $C_p$ & $\times$ \\
\hline\hline
\end{tabular}, \ \ \ 
\begin{tabular}{||c||c||}
\hline\hline $*$ & ${\cal C}({\cal M})$\\ \hline\hline ${\cal C}({\cal
M})$ & $\times$ \\ \hline\hline
\end{tabular}  \ \ \right)$$
}

\vskip -0.3cm
$$
\begin{tabular}{||c||c||c||c||c||}
\hline\hline $***$ & $C_p^1$ & $C_p^2$ & $C_p^3$ & $C_p^4$ \\
\hline\hline $C_p^1$ & $\Leftarrow \ \ \Rightarrow$ & $\nLeftarrow \ \
\Rightarrow$ & $\nLeftarrow \Rightarrow$ & $\nLeftarrow \ \
\Rightarrow$ \\ \hline\hline $C_p^2$ &$\Leftarrow \ \ \nRightarrow$
&$\Leftarrow \ \ \Rightarrow$ &$\nLeftarrow \ \ \nRightarrow$
&$\nLeftarrow \ \ \Rightarrow$ \\ \hline\hline $C_p^3$ &$\Leftarrow \
\ \nRightarrow$ &$\nLeftarrow \ \ \nRightarrow$ &$\Leftarrow \ \
\Rightarrow$ &$\nLeftarrow \ \ \Rightarrow$ \\ \hline\hline $C_p^4$
&$\Leftarrow \ \ \nRightarrow$ &$\Leftarrow \ \ \nRightarrow$
&$\Leftarrow \ \ \nRightarrow$ & $\Leftarrow \ \ \Rightarrow$\\
\hline\hline
\end{tabular}
$$

\smallskip
$$
\begin{tabular}{||c||c||c||c||c||}
\hline\hline $***$ & $C^1(M)$ & $C^2(M)$ & $C^3(M)$ & $C^4(M)$ \\
\hline\hline $C_p^1$ &$\nLeftarrow \ \ \Rightarrow$ &$\nLeftarrow \ \
\Rightarrow$ &$\nLeftarrow \ \ \Rightarrow$ & $\nLeftarrow \ \
\Rightarrow$\\ \hline\hline $C_p^2$ &$\nLeftarrow \ \ \nRightarrow$
&$\nLeftarrow \ \ \nRightarrow$\footnotemark[11] &$\nLeftarrow \ \
\nRightarrow$\footnotemark[11] &$\nLeftarrow
\nRightarrow$\footnotemark[11] \\ \hline\hline $C_p^3$ &$\nLeftarrow \
\ \nRightarrow$ & $\nLeftarrow \ \ \nRightarrow$\footnotemark[11] &
$\nLeftarrow \ \ \nRightarrow$\footnotemark[11] & $\nLeftarrow \ \
\nRightarrow$\footnotemark[12] \\ \hline\hline $C_p^4$ &$\nLeftarrow \
\ \nRightarrow$ &$\nLeftarrow \ \ \nRightarrow$\footnotemark[11]
&$\nLeftarrow \ \ \nRightarrow$\footnotemark[11] &$\nLeftarrow \ \
\nRightarrow$\footnotemark[11] \\ \hline\hline
\end{tabular}
$$

\smallskip
$$
\begin{tabular}{||c||c||c||c||c||}
\hline\hline $***$ & ${\cal C}^1({\cal M})$ & ${\cal C}^2({\cal M})$ &
${\cal C}^3({\cal M})$ & ${\cal C}^4({\cal M})$ \\ \hline\hline ${\cal
C}^1({\cal M})$ &$\Leftarrow \ \ \Rightarrow$
&$\nLeftarrow\footnotemark[12] \Rightarrow$
&$\nLeftarrow\footnotemark[12] \Rightarrow$ &
$\nLeftarrow\footnotemark[12] \Rightarrow$\\ \hline\hline ${\cal
C}^2({\cal M})$ &$\Leftarrow \ \ \nRightarrow$\footnotemark[12]
&$\Leftarrow \ \ \Rightarrow$ & $\nLeftarrow\footnotemark[13]
\nRightarrow$\footnotemark[14]&$\nLeftarrow\footnotemark[13]
\Rightarrow$ \\ \hline\hline ${\cal C}^3({\cal M})$ &$\Leftarrow \ \
\nRightarrow\footnotemark[12]$ &$\nLeftarrow\footnotemark[14]
\nRightarrow\footnotemark[13]$ & $\Leftarrow \ \
\Rightarrow$&$\nLeftarrow\footnotemark[14] \Rightarrow$ \\
\hline\hline ${\cal C}^4({\cal M}))$ &$\Leftarrow \ \
\nRightarrow\footnotemark[12]$ &$\Leftarrow \ \
\nRightarrow\footnotemark[13]$ &$\Leftarrow \ \
\nRightarrow\footnotemark[14]$ &$\Leftarrow \ \ \Rightarrow$ \\
\hline\hline
\end{tabular}
$$

\bigskip

Our examples are intended to explain only the main nontrivial (non)
implication links above. Those not in the articles are easier to find.

\footnotetext[11]{Example 13 shows all: $C_p^j \nRightarrow C^j$,
$2\leq j\leq 4$.}  \footnotetext[12]{Example 8 shows all: ${\cal
C}^j({\cal M}) \nRightarrow {\cal C}^1({\cal M})$, $2\leq j\leq 4$.}
\footnotetext[13]{Example 11.}  \footnotetext[14]{Example 12.}

\smallskip
\noindent
{\bf k. Nonanalytic behavior of $(z,\bar{w})\mapsto \X_{z,\bar{w}}'$.}

Two examples are given to exhibit it.

\smallskip
\noindent
{\bf Example 13.} There exist $f$, $M$, $M'$ with $f$ nonalgebraic
such that the function ${\cal M} \ni (z,\bar{w}) \mapsto {\rm
dim}_{f(z)}\X_{z,\bar{w}}' \in \N$ is not upper semi-continuous at
$0$. Explicitely, take $M: z_4=\bar{z}_4+iz_1\bar{z}_1$ in
$\C^2_{(z_1, z_4)}$, \begin{equation} M': \ \ \
z_4'=\bar{z}_4'+iz_1'\bar{z}_1'+ iz_3' \bar{z}_2'+i\bar{z}_3' z_2'+
i{z'}_1^2\bar{z}_3' \bar{z}_2' +i \bar{z'}_1^2 z_3' z_2'
\end{equation} and $f(z_1,z_4)=(z_1, z_4\sin^3 z_1, 0, z_4)$.  Here,
${\rm dim}_0 \X_{0,0}' =0$, ${\rm dim}_{f(z)} \X_{z,\bar{w}}' =1$
$\forall \ z \neq 0$, $z\in Q_{\bar{w}}$. (Idem for $\M_{z,\bar{w}}'$
instead.)  See Section 4.

\smallskip
This example therefore shows that $\X_{z,\bar{w}}'$ cannot be written
as $\{z' \ \! {\bf :} \ \!  \lambda(z,\bar{w}, z')=0\}$ for
holomorphic $\lambda\in {\cal V}_{\C^n}(0) \times {\cal V}_{\C^n}(0)
\times {\cal V}_{\C^{n'}}(0)$ in general.

\smallskip
\noindent
{\bf Example 14.} There exist $f$, $M$, $M'$ all algebraic such that
if $\Sigma$ denotes the set of $(z,\bar{w}, z_1) \in {\cal M} \sharp
{\cal M}$ in a neighborhood of which $(*)$ of Proposition 5 is not
satisfied, then $\Sigma$ is not a complex analytic subset of ${\cal M}
\sharp {\cal M}$ but a real analytic set.

\smallskip
\noindent
{\bf l. Globalization of $r_{M'}$, $r_{M'}^2$.}

First, we notice that $r_M(E) (= r_M^U(E)$, localized in $U$, {\it cf}
[Z98]), could have been defined globally as \begin{equation}
r_M^{\C^n}(E)=\{w\in \C^n\ \! {\bf :} \ \!  Q_{\bar{w}} \supset E \}
\end{equation} because $\rho, \rho'$ are polynomials so that it holds
that

\smallskip
\noindent
{\bf Theorem 1''.} {\it If ${\rm dim}_{f(z)}[r_{M'}^{\C^{n'}}
(f(Q_{\bar{z}})) \cap (r_{M'}^{\C^{n'}})^2(f(Q_{\bar{w}}))]=0$
$\forall z, w \in {\cal V}_{\C^n}(0)$, $z\in Q_{\bar{w}}$, then $f$ is
algebraic.}

\smallskip
\noindent
(identical proof).

Then Theorem 1'' can be more general than Theorem 1 because

\smallskip
\noindent
{\bf Example 15.} There exist $f$, $M$, $M'$, $U$, $U'$ such that

1. ${\rm dim}_{f(z)}[r_{M'}^{U'} (f(Q_{\bar{z}})) \cap
(r_{M'}^{U'})^2(f(Q_{\bar{w}}))]\geq 1$ $\forall z, w \in U$, $z\in
Q_{\bar{w}}$ and

2. ${\rm dim}_{f(z)}[r_{M'}^{\C^{n'}} (f(Q_{\bar{z}})) \cap
(r_{M'}^{\C^{n'}})^2(f(Q_{\bar{w}}))]=0$ $\forall z, w \in U$, $z\in
Q_{\bar{w}}$.

\noindent
Explicitely (see Section 4), take $U=\Delta^2$, $U'=\Delta^4$, the
hypersurface $M: z_4=\bar{z}_4 +iz_1\bar{z}_1$, $f(z_1, z_4)=(z_1,0,0,
z_4)\in \C^4$, and \begin{equation} M':
z_4'=\bar{z}_4'+iz_1'\bar{z}_1' +i{z'}_1^2(1+\bar{z}_3')\bar{z}_3'+
i\bar{z'}_1^2(1+z_3')z_3'+ iz_2' z_3' \bar{z'}_2^2+i \bar{z}_2'
\bar{z}_3' {z'}_2^2.  \end{equation} Conversely, Theorem 1 can be more
general than Theorem 1'' (exercise left to the reader).

\smallskip
\noindent
{\bf m. About $\Z_{z,\bar{w}}'$, $\M_{z,\bar{w}}'$.}

It is not difficult to see that all the positive results Theorems 1,
1', 1'' concerning $\X_{z,\bar{w}}'$ extend to be satisfied by
$\Z_{z,\bar{w}}'$ and by $\M_{z,\bar{w}}'$ also once one has
established the following result analogous to Proposition 5 and with
even a better property: the bad set ${\cal N}$ below is a complex
analytic subvariety of ${\cal M}$.

\smallskip
\noindent
{\bf Proposition 16.} {\it There exists a standard constructive way of
finding a variety $\widetilde{\S}^1=\{(z,\bar{w}, z')$ ${\bf :} \ \!
\tilde{r}(z,\bar{w}, z')= 0\}$ contained in $\S^1$ with
$\tilde{r}_j(z,\bar{w}, z') \in {\cal A}_n (U) \times \overline{\cal
O}_n(U)\times {\cal A}_{n'}(U')$, $1 \leq j \leq \widetilde{J} $, $
\widetilde{J} \geq J$, such that \begin{equation} \Gamma r(f) =
\{(z,\bar{w}, f(z)) \ \! {\bf :} \ \!  (z,\bar{w}) \in {\cal M} \}
\subset \widetilde{\S}^1_{{\cal M}}, \end{equation}
$$\widetilde{\S}^1_{{\cal M}} = \{(z,\bar{w}, z') \ \! {\bf :} \ \!
\tilde{r}(z,\bar{w}, z') = 0 , (z,\bar{w}) \in {\cal M}\}
$$
and such that there exist a Zariski open subset ${\cal D}_{\cal M}:=
{\cal M} \backslash {\cal N}$ of ${\cal M}$, ${\cal N}\subset {\cal
M}$ complex analytic, ${\rm dim}_{\C} {\cal N} \leq 2m + d - 1$, and
an integer $n_1'$, $0\leq n_1'\leq n'$ such that

$(*)$ $\forall \ p \in {\cal D}_{\cal M}$, \ $\exists\ {\cal U}_p =
{\cal V}_{\cal M}(p)$, $\exists U_{p'}' = {\cal
V}_{\C^{n'}}(p'=f(p))$, such that $\widetilde{\S}_{{\cal U}_p}
^1:=({\cal U}_p \times U_{p'}') \cap \widetilde{\S}^1$ is smooth and
the projection $\pi': \widetilde{\S}^1_{{\cal U}_p}\to U'$ is of
constant rank $n_1'$.

Consequently, $(*)$ implies that

$(**)$ $\forall \ p \in D_M := {\cal D}_{\cal M} \cap M$, $\exists \
U_p = {\cal V}_{\C^n}(p)$, $\exists U_{p'}' = {\cal
V}_{\C^{n'}}(p'=f(p))$, $\exists K \in \N_*$, $\exists (s_k)_{1\leq
k\leq K}$, $s_k \in \overline{\cal O}_n({U}_p) \times {\cal O}_n (U_p)
\times {\cal A}_{n'}(U_{p'}')$ such that $\forall \ z_1 \in
Q_{\bar{w}}$ if $\W_{w,\bar{z}_1}':=
\widetilde{\S}^1_{w,\bar{z}_1}\cap (U_p \times \overline{U}_p \times
U_{p'}') $ for any $(w,\bar{z}_1)\in {\cal M}$, and if
$r_{M'}^{U_{p'}'}(E'):= \{w' \in U_{p'}' \ \! {\bf :} \ \!
Q_{\bar{w}'}' \supset E'\}$, then \begin{equation} r_{M'}^{U_{p'}'}
(\W_{w,\bar{z}_1}) = \{(\bar{w}, z_1, z') \in \overline{U}_p \times
U_p \times U_{p'}'\ \! {\bf :} \ \!  s(\bar{w}, z_1, z') =0\}
\end{equation}
}

\vskip -0.2cm Proposition 16 shows that after shrinking the first
reflection variety $\S^1$ to $\widetilde{\S}^1$, one is enabled to
satisfy property $(*)$ above ({\it cf.} ``holomorphic families'' in
[Z98]) of which the computation of the second reflection
$r_{M'}(\W_{w,\bar{z_1}}')$ {\it after localisation} in a smaller open
subset $U_{p'}'$ is {\it easier} than the computation of
$r_{M'}(\S_{w,\bar{z}_1}^1)$ and gives its desired analytic dependence
with respect to the parameters $(w,\bar{z}_1)$ as written in $(**)$
above eq. (1.37).

We recall that our examples show that there is a serious difference
between Proposition 5 and Proposition 16 and a serious difference
between applying operators $r_{M'}^{U_{p'}'}$ or $r_{M'}^{U'}$ or
$r_{M'}^{\C^{n'}}$.

Finally, we obtain from Proposition 16 the Theorems 1 and 1' with
$\Z_{z,\bar{w}}'$ and with $\M_{z,\bar{w}}'$ instead of
$\X_{z,\bar{w}}'$. \hfill $\square$

\smallskip
\noindent
{\bf Acknoweledgement.} The author would like to thank S. Damour for
fruitful conversations.

\section{Proof of Theorem 1}

\def\theequation{2.\arabic{equation}} \setcounter{equation}{1}

\noindent
{\bf Lemma 2.1.} {\it $(i)$ For any set $E\subset U$, $E\cap r_M(E)
\subset M$ and $E\subset r_M(r_M(E))$.

$(ii)$ $z\in Q_{\bar{w}}$ iff $w\in Q_{\bar{z}}$, $z\in Q_{\bar{z}}$
iff $z\in M$, $\{w\in U\ \! {\bf :} \ \! Q_{\bar{w}} \supset
E\}=\cap_{w\in E} Q_{\bar{w}}$.

$(iii)$ $f(Q_{\bar{z}}) \subset Q_{\overline{f(z)}}'$.

$(iv)$ $\rho(z,\bar{w})=0$ iff $\rho(w,\bar{z})=0$.}

\smallskip
{\it Proof.}  (ii), (iii) and (iv) are classical. Prove (i). If $e\in
E$ and $e\in r_M(E)= \cap_{w\in E} Q_{\bar{w}}$ then $e\in
Q_{\bar{e}}$, so $e\in M$ by (ii), {\it i.e.} $E\cap r_M(E) \subset
M$.  Furthermore, by construction of $r_M(E)$, \begin{equation}
r_M(r_M(E))=\cap \{Q_{\bar{z}} \ \! {\bf :} \ \!  z\in r_M(E) \} =\cap
\{Q_{\bar{z}} \ \! {\bf :} \ \!  Q_{\bar{z}} \supset E\} \supset E. \
\ \square \end{equation}

Let ${\cal M}=\{(z,\bar{w})\in U\times U\ \! {\bf :} \ \!
\rho(z,\bar{w})=0\}$. Let $\underline{\cal L}_l=\sum_{j=1}^n
a_{j,l}(z,\bar{w})\frac{\partial }{\partial \bar{w}_j}$, $1\leq l \leq
m$, be tangent vectors to ${\cal M}$ which are the complexifications
of a basis of tangent vectors $\overline{L}_l= \sum_{j=1}^{n}
a_{j,l}(z,\bar{z}) \frac{\partial }{\partial \bar{z}_j}$, $1\leq l\leq
m$ generating $T^{0,1} M$ with polynomial coefficients and which
commute. Let ${\cal O}_n$ (resp. ${\cal A}_n$) denote the ring of
holomorphic functions in ${\cal V}_{\C^n}(0)$ (resp.  {\it holomorphic
polynomials}). Here, ${\cal V}_{\C^n}(0)=U$ after possible shrinking.

\def\theequation{2.\arabic{equation}} \setcounter{equation}{3}

\smallskip
\noindent
{\bf Lemma 2.3.} {\it There exist $J\in \N_*$ and functions
$r_j(z,\bar{w},z')\in {\cal A}_n \times \overline{\cal O}_n \times
{\cal A}_{n'}$, $1\leq j\leq J$, such that $\forall \ (z,\bar{w}) \in
{\cal M}$ \begin{equation} r_{M'} (f(Q_{\bar{z}}))=\{z'\in U' \ \!
{\bf :} \ \!  r_j(z,\bar{w}, z')=0, 1\leq j\leq J\}.  \end{equation}}

{\it Remark.} Sets $\{r_j(z,\bar{w}_1, z')=0\}$ and
$\{r_j(z,\bar{w}_2, z')=0\}$ for different $w_1, w_2$ such that $z\in
Q_{\bar{w}_1}$, $z\in Q_{\bar{w}_2}$, coincide and are equal to
$r_{M'} (f(Q_{\bar{z}}))$.

\smallskip
{\it Proof.} By definition, $r_{M'}(f(Q_{\bar{z}}))=\{w'\in U'\ \!
{\bf :} \ \!  \rho'(f(w), \bar{w}')=0, \forall \ w\in Q_{\bar{z}}\}$.
Using Lemma 2.1 (iv), $r_{M'}(f(Q_{\bar{z}}))=\{z'\in U' \ \! {\bf :}
\ \!  \rho'(z', \bar{f}(\bar{w}))=0, \forall \ w\in
Q_{\bar{z}}\}$. Equivalently, $Q_{\bar{z}} \ni w \mapsto \rho'(z',
\bar{f}(\bar{w}))\in \C^{d'}$ vanishes identically as an
antiholomorphic map of $w$ defined on the complex algebraic manifold
$Q_{\bar{z}}$. Because of identity principle, this is equivalent to
$\underline{\cal L}^{\gamma} (\rho'(z',\bar{f}(\bar{w})))=0$ $\forall
\ \gamma\in \N^m$. Put $r_{\gamma}(z,\bar{w}, z'):= \underline{\cal
L}^{\gamma} (\rho'(z', \bar{f}(\bar{w})))$. Then $r_{\gamma} \in {\cal
A}_n\times \overline{\cal O}_n \times {\cal A}_{n'}'$.  By
noetherianity, a finite number $J$ of $r_{\gamma}$'s defines
$r_{M'}(f(Q_{\bar{z}}))$.  \hfill $\square$

\def\theequation{2.\arabic{equation}} \setcounter{equation}{5}

\smallskip
\noindent
{\bf Lemma 2.5.} {\it There exist $K\in \N_*$ and polynomials
$s_k(z')$ (depending on $\bar{w}$), $1\leq k\leq K$, such that
$r_{M'}^2(f(Q_{\bar{w}}))=\{z'\in U' \ \! {\bf :} \ \!  s_k(z')=0,
1\leq k\leq K\}$.}

\smallskip
{\it Proof.} Simply because $r_{M'}^2(f(Q_{\bar{w}}))$ is algebraic,
by (1.1). \hfill $\square$

\smallskip
{\it End of proof of Theorem 1.} Fix $\bar{w}$. We prove that
$f|_{Q_{\bar{w}}\cap V}$ is algebraic. Indeed

1. $\forall \ z\in Q_{\bar{w}}$, $r_j(z,\bar{w}, f(z))\equiv 0$,
$1\leq j\leq J$ and $s_k(f(z))\equiv 0$, $1\leq k\leq K$.

2. The set $\Phi:=\{z' \ \! {\bf :} \ \!  r_j(z,\bar{w}, z') =0, 1\leq
j\leq J \ \hbox{and} \ s_k(z')=0, 1\leq k\leq K\}$ is zero dimensional
at each point $f(z)$, $z\in V$, $z\in Q_{\bar{w}}$.

\smallskip
By Theorem 5.3.9 in [BER99] (in the algebraic case), there exist
Weierstrass polynomials $P_j(z, z_j') ={z_j'}^{N_j} +\sum_{1\leq k\leq
N_j} A_{k,j}(z) {z_j'}^{N_j-k}$, $A_{k,j} \in \C[z]$, $1\leq j\leq
n'$, such that $\Phi$ is contained in $\Psi=\{(z,z') \in V\times V' :
P_j(z, z_j')=0, 1\leq j\leq n'\}$. As $(z,f(z))\in \Phi$, we obtain
that \begin{equation} f_j(z)^{N_j} +\sum_{1\leq k\leq N_j} A_{k,j}(z)
f_j^{N_j-k}(z)\equiv 0, \ \ \ \ \ z\in Q_{\bar{w}} \cap V, \ 1\leq
j\leq n'.  \end{equation} Equation above yields at once that each map
$Q_{\bar{w}}\cap V \ni z\mapsto f_j(z) \in \C$ is holomorphic {\it
algebraic}, $1\leq j\leq n'$. To conclude, apply Theorem 2.  \hfill
$\square$

\section{Proof of Proposition 5}

Proposition 5 relies upon the following statement (denseness of $D_M$
is then clear and $(*)\Rightarrow (**)$ also):

\smallskip
\noindent
 {\bf Proposition 3.1.} {\it If ${\rm dim}_{f(z)} \X_{z,\bar{z}}'=0$
$\forall \ z\in M\cap {\cal V}_{\C^n}(0)=M\cap V$, then there exists
$p\in M\cap V$ arbitrarily close to $0$ such that ${\rm
dim}_{f(z)}\X_{z,\bar{w}}'=0$, $\forall \ z, w \in {\cal
V}_{\C^n}(p)$, $z\in Q_{\bar{w}}$ and $(*)$ of Proposition 5 holds in
$U_p ={\cal V}_{\C^n}(p)\subset V$.}

\smallskip
As we shall see, the main difficulty is that there does not
necessarily exist holomorphic equations $\lambda(z,\bar{w}, z')$ such
that $\X_{z,\bar{w}}'=\{z'\in U' \ \! {\bf :} \ \!  \lambda(z,\bar{w},
z')=0\}$ as for example like for $r_{M'}(f(Q_{\bar{z}}))=\{z'\in U' \
\! {\bf :} \ \!  r(z,\bar{z}, z')=0\}$, $z\in M$, $z\in Q_{\bar{z}}$,
equations which would readily imply the desired conclusion by
semi-continuity of the fiber dimension of a holomorphically
parametrized family of complex analytic sets.  To get such a local
parametrized family, we shall have to shift $p$ from a certain number
of images by holomorphic maps of complex analytic sets by keeping
$p\in M$ in the same time and our proof gives that the set of $p\in M$
in a neighborhood of which $\X_{z,\bar{w}}'$ should be holomorphically
parametrized is a dense open subset of $M$. It will also clearly show
that the bad set can be at least as worst as a subanalytic set.

We will work out Proposition 3.1 with $M$, $M'$ of class ${\cal
C}^{\omega}$.

For that purpose, let $\V_z':=r_{M'} (f(Q_{\bar{z}}))=\{z'\in U' \ \!
{\bf :} \ \!  r(z,\bar{w}, z')=0\}$ (in vectorial notations,
$r=(r_1,\ldots,r_J)$) so that
$r_{M'}^2(f(Q_{\bar{w}}))=r_{M'}(\V_w')$. Let us recall that the
representation of $\V_z'$ by holomorphic equations $r(z,\bar{w}, z')$
gives the same set $\V_z'$ for any choice of $(z,\bar{w})\in {\cal M}$
({\it cf.} Lemma 2.3). Therefore the introduction of a third point
$z_1\in Q_{\bar{w}}$ yields a representation $\V_w'=\{z'\in U'\ \!
{\bf :} \ \!  r(w,\bar{z}_1, z')=0\}$.

From now on, we let $z, w, z_1\in U$, $z\in Q_{\bar{w}}$, $w\in
Q_{\bar{z}_1}$, and we denote $\S_{z,\bar{w}}^1$, $\S_{w,\bar{z}_1}^1$
instead of $\V_z'$, $\V_w'$. This is justified by the fact that
although the set $\{z'\in U'\ \! {\bf :} \ \!  r(z,\bar{w}, z')=0\}$
does not depend on $\bar{w}$, the {\it equations}
$r_{\gamma}(z,\bar{w}, z')= \underline{\cal L}^{\gamma} \rho'(z',
\bar{f}(\bar{w}))=0$ {\it do really depend} on $\bar{w}$.  (Inspect
for instance the identity map $\C^2 \to \C^2$, $M\to M$,
$M=\{z_2=\bar{z}_2+iz_1\bar{z}_1\}$.)  The notation $\S_{z,\bar{w}}^1$
simply means a section over $(z,\bar{w})$ of the set $\S^1$ even if
$(z,\bar{w}) \not\in {\cal M}$.

\def\theequation{3.\arabic{equation}} \setcounter{equation}{1} We then
write \begin{equation} r_{M'}(\S_{w,\bar{z}}^1)= \{z'\in U' \ \! {\bf
:} \ \! \rho'(z', \overline{\zeta}')=0 \ \forall \ \zeta' \ {\rm
satisfying} \ r(w,\bar{z}_1, \zeta')=0\}.  \end{equation} We shall
establish that there exist points $p=(z_p, \bar{z}_p) \in M$
arbitrarily close to $0$, neighborhoods $U_p = {\cal V}_{\C^n}(p)$,
$U_{p'}'= {\cal V}_{\C^{n'}}(p'=f(p))$ and holomorphic equations
$s(\bar{w}, z_1, z')$ near $(z_p, \bar{z}_p, f(z_p))$ in $U_p \times
\overline{U}_p \times U_{p'}'$ such that

1. $r_{M'}(\S_{w,\bar{z}_1}^1) \subset \{z'\in U_{p'}' \ \! {\bf :} \
\!  s(\bar{w}, z_1, z')=0\}$ $\forall \ w, z_1 \in U_p$, $w\in
Q_{\bar{z}_1}$ and

2. $r_{M'}(\S_{w,\bar{w}}^1)=\{z'\in U_{p'}' \ \! {\bf :} \ \!
s(\bar{w}, w, z')=0\}$ $\forall w\in U_p \cap M$

\smallskip\def\theequation{3.\arabic{equation}}
\setcounter{equation}{3}
\noindent
 {\bf Lemma 3.3.} {\it If the above conditions 1-2 are fulfilled, then
${\rm dim}_{f(z)} \X_{z,\bar{w}}'=0$ $\forall \ (z,\bar{w})\in {\cal
V}_{{\cal M}}(p)$.}

\smallskip
{\it Proof.} Indeed, let $\Y_{z,\bar{w}, z_1}'=\{z'\in U' \ \! {\bf :}
\ \!  r(z,\bar{w}, z')=0 , s(\bar{w}, z_1, z')=0\}\supset
\X_{z,\bar{w}}' \ni f(z)$. By assumption, ${\rm dim}_{f(p)}
\Y_{p,\bar{p}, p}' =0$.  Since $\Y_{z,\bar{w}, z_1}'$ is
holomorphically parametrized, then ${\rm dim}_{f(z)} \Y_{z,\bar{w}
,z_1}'=0$, $\forall \ z, w, \bar{z}_1$ in some small neighborhood of
$p$ in $\C^n$ with $z\in Q_{\bar{w}}$, $w\in Q_{\bar{z}_1}$.
Therefore ${\rm dim}_{f(z)} \X_{z,\bar{w}}'=0$ also.  \hfill $\square$

\smallskip
We have obtained a dense open set $D_M$ where the above conditions 1
and 2 are fulfilled.  To complete $(*)$ of Proposition 5, it suffices
to take $\S^2=$ the irreducible component of $\Y'=\{(z,\bar{w},
z_1,z')\in (({\cal M} \sharp {\cal M} ) \cap (U_p \times
\overline{U}_p \times U_p))\times U_{p'}' \ \! {\bf :} \ \!
r(z,\bar{w}, z')=0, s(\bar{w}, z_1, z')=0\}$ containing the graph of
$f$ over $({\cal M}\sharp {\cal M})\cap (U_p\times
\overline{U}_p\times U_p)$. $\S^2$ is defined by similar partially
polynomial equations $r\in {\cal A}_n(U_p)\times \overline{\cal
O}(U_p) \times {\cal A}_{n'} (U_{p'}')$, $s\in \overline{\cal
O}_n(U_p) \times {\cal O}_n(U_p)\times {\cal A}_{n'}(U_{p'}')$.  The
graph of $f$ is in fact a local irreducible component of $\S^2$, for
reasons of dimension. It is well-known that in the algebraic case,
local irreducible components cannot be defined by polynomial equations
at every point. Therefore, the component is a smooth complex manifold
defined by the same partially polynomial equations equations
$r(z,\bar{w}, z')=0$ and $s(\bar{w}, z_1, z')=0$ near each point
belonging to some Zariski open subset of ${\cal M}\sharp {\cal
M}$. This Zariski open set intersects $M\sharp M\cong M$ in a Zariski
open subset of $M\cap U_p$, because the $r$ depends only on $z,
\bar{w}$ (not on $z_1$) and $s$ depends only on $\bar{w}, z_1$ (not on
$z$) and because $M$ is maximally real in ${\cal M}$.  Therefore,
after perharps shifting the points in $U_p\cap M$ and thus shrinking
$D_M$, $(*)$ of Proposition 5 is satisfied.  \hfill $\square$

\smallskip

Consequently, we are reduced to a statement which we now formalize in
an independent fashion leaving Segre varieties, see Lemma 3.7.

Let $\Delta$ be the unit disc in $\C$.

Let $\kappa\in \N_*$, $n\in \N_*$, $\nu_1\in \N_*$, $g:
\Delta^{\kappa} \times \Delta^n\to \C^{\nu_1}$ be a holomorphic power
series mapping converging uniformly in $(2\Delta)^{\kappa+n}$,
$(t,z)\mapsto g(t,z)$, let \begin{equation} F=\{(t,z)\in
\Delta^{\kappa} \times \Delta^n \ \! {\bf :} \ \!  g(t,z) =0\}.
\end{equation} Assume that there exists $\lambda: \Delta^{\kappa} \to
\Delta^n$ holomorphic, converging in $(2\Delta)^{\kappa}$ such that
$(t, \lambda(t))\in F$, $\forall \ t\in \Delta^{\kappa}$, hence
$\pi(F)= \Delta^{\kappa}$, where $\pi: \Delta^{\kappa} \times \Delta^n
\to \Delta^{\kappa}$, $(t, z)\mapsto t$. Let $n'\in \N_*$, $\nu_2\in
\N_*$, let $\rho: \Delta^{n'}\times \Delta^n \to \C^{\nu_2}$ be a
holomorphic series converging in $(2\Delta)^{n'+n}$ and denote for
$t\in \Delta^{\kappa}$ \begin{equation} G_F[t]=\{z'\in \Delta^{n'} \
\! {\bf :} \ \!  \rho(z', z)= 0 \ \forall \ z \ \ s.t. \ \ g(t,z)=0\}.
\end{equation} Let $I^{\kappa}\subset \Delta^{\kappa}$ be the
maximally real set $I^{\kappa}=(-1,1)^{\kappa}$, $I=(-1,1)$.

Introduce the complex filtration $F=F_1\supset F_2 \supset F_3 \supset
\cdots \supset F_{a+1}=\emptyset$, $F_a \neq \emptyset$, $a\geq 1$, of
$F$ by singular subspaces: $F_{i+1}=F_{i, sing}$. Assume also that
\begin{equation} F_{\alpha}=\{(t,z)\in \Delta^{\kappa}\times \Delta^n
\ \! {\bf :} \ \! g_{\alpha}(t,z)=0\}, \ \ \ \ \ 1\leq \alpha\leq a,
\end{equation} with $g_{\alpha}: (2\Delta)^{\kappa+n} \to
\C^{\nu_{1\alpha}}$ like $g$ and that all irreducible components of
$F_{\alpha}$ are defined analogously.  We denote by
$\Delta_{\kappa}(\underline{t}, \underline{\varepsilon})$,
$\underline{\varepsilon} >0$, the polydisc of center $\underline{t}$,
radius $\underline{\varepsilon}$, with $\underline{t} \in
\Delta^{\kappa}$, $\underline{\varepsilon} << {\rm dist}
(\underline{t}, b\Delta^{\kappa})$.  It remains to establish:

\smallskip\def\theequation{3.\arabic{equation}}
\setcounter{equation}{7}
\smallskip
\noindent
 {\bf Lemma 3.7.} {\it There exist $\underline{t}\in I^{\kappa}$,
$\underline{\varepsilon} >0$ and holomorphic equations $s(t,z')$ in
$\Delta_{\kappa}(\underline{t}, \underline{\varepsilon})\times
\Delta^{n'}$ such that

1. $G_F[t] \subset \{z'\in \Delta^{n'} \ \! {\bf :} \ \! s(t,z')=0\}$
\ \ \ $\forall\ t\in \Delta_{\kappa}(\underline{t},
\underline{\varepsilon})$ and

2. $G_F[t] =\{z'\in \Delta^{n'} \ \! {\bf :} \ \!  s(t,z') =0\}$ \ \ \
$\forall \ t\in \Delta_{\kappa}(\underline{t},
\underline{\varepsilon}) \cap I^{\kappa}$.}

\smallskip
{\it Proof.} First, $G_F[t] = G_{F_{1, reg}}[t] \cap G_{F_{2,reg}}[t]
 \cap \cdots \cap G_{F_{a,reg}}[t]$ $(F_a=F_{a,reg})$. Also, $\forall
 \ \alpha$, $1\leq \alpha \leq a$, $G_F[t]= G_{F_{1,reg}}[t]\cap
 \cdots \cap G_{F_{\alpha-1, reg}}[t] \cap G_{F_{\alpha}}[t]$. Let
 $F_{\alpha}=\cup_{\beta=1}^{b_{\alpha}} F_{\alpha \beta}$,
 $b_{\alpha}\in \N_*$, denote the decomposition of $F_{\alpha}$ into
 irreducible components. Then also \begin{equation} G_F[t] =
 \bigcap_{1\leq \alpha \leq a, 1\leq \beta \leq b_{\alpha}}
 G_{F_{\alpha\beta}}[t].  \end{equation} We have (see [C], Chapter 1):

\hspace{1cm} $F_{\alpha\beta, reg} \backslash (\cup_{\gamma \neq
\beta} F_{\alpha \gamma}) \subset F_{\alpha, reg}$,

\hspace{1cm} $F_{\alpha\beta, reg} \cap (\cup_{\gamma\neq \beta}
F_{\alpha\gamma})$ is a proper analytic subset of $F_{\alpha\beta,
reg}$,

\hspace{1cm} $F_{\alpha, reg}= \cup_{1\leq \beta \leq b_{\alpha}}
(F_{\alpha\beta, reg} \backslash (\cup_{\gamma \neq \beta} F_{\alpha
\gamma}))$ and

\hspace{1cm} $G_{F_{\alpha, reg}}[t] \cap G_{F_{\alpha+1}}[t] \subset
\cap_{1\leq \beta \leq b_{\alpha}} G_{F_{\alpha\beta, reg}}[t]$.

\noindent
This yields \begin{equation} G_F[t]=\bigcap_{1\leq \alpha \leq a,
                      1\leq \beta \leq b_{\alpha}}
                      G_{F_{\alpha\beta,reg}}[t].  \end{equation}
                      Denote now by $F_1, \ldots, F_c$,
                      $c=b_1+\cdots+b_a\in \N_*$ the
                      $F_{\alpha\beta}$'s which are irreducible. So
                      $G_F[t] = \cap_{1\leq \gamma \leq c}
                      G_{F_{\gamma,reg}}[t]$.

Let $F$ be one of the $F_{\gamma}$'s, $1\leq \gamma \leq c$.

Now, we come to a dichotomy.

Either the generic rank satisfies \begin{equation} {\rm gen \ rk}_{\C}
                      (\pi|_{F_{reg}}) = \kappa \ \ \ \ \ {\rm or} \ \
                      \ \ \ {\rm gen \ rk}_{\C} (\pi|_{F_{reg}}) <
                      \kappa.  \end{equation}

\smallskip\def\theequation{3.\arabic{equation}}
\setcounter{equation}{11}
\smallskip
\noindent
{\bf Lemma 3.11.} {\it Let $F=$ one of the $F_{\gamma}$. If $ {\rm gen
\ rk}_{\C} (\pi|_{F_{reg}}) < \kappa$, then the {\rm closed} set
$\overline{\pi}(F):=\pi(\overline{F} \subset
\overline{\Delta}^{\kappa} \times \overline{\Delta}^{n}) \subset
\overline{\Delta}^{\kappa}$ (here, $\overline{\ \ }$ denotes closure)
is contained in a countable union $\cup_{\nu\in \N_*} A_{\nu}$ of
analytic sets $A_{\nu} \cong\Delta^{\lambda_{\nu}}$ with $0\leq
\lambda_{\nu} < \kappa$, $\nu\in \N_*$.}

\smallskip
{\it Proof.}  Let $F$ be a $F_{\gamma}$ with ${\rm gen \ rk}_{\C}
(\pi|_{F_{reg}}) < \kappa$.  Since $F$ is defined over
$(2\Delta)^{\kappa+n}$ and irreducible, paragraph 3.8 in [C]
applies. \hfill $\square$

\smallskip
Hence the Lebesgue measure
$\lambda_{2\kappa}(\overline{\pi}(F))=0$. Furthermore, $\forall \
\nu$, $\lambda_{\kappa} (\overline{I}^{\kappa} \cap A_{\nu}) = 0$,
since $I^{\kappa}$ is maximally real. Hence
$\lambda_{\kappa}(\overline{\pi}(F)\cap I^{\kappa})=0$.  Thus there
exists an open dense subset $B_F$ of $I^{\kappa}$ such that $\forall \
\underline{t} \in B_F$, $\exists$ an open neighborhood ${\cal
V}_{\Delta^{\kappa}} (\underline{t})$ with ${\cal V}_{\Delta^{\kappa}}
(\underline{t}) \cap \overline{\pi}(F) = \emptyset$.

{\it Consequently, all irreducible components $F_{\gamma}$ such that
${\rm gen \ rk}_{\C} (\pi|_{F_{\gamma,reg}}) < \kappa$ can be
forgotten.} Indeed, for almost all $t\in \Delta^{\kappa}$,
$F_{\gamma}[t]=\emptyset$, so for such $t$, $F_{\gamma}$ makes no
contribution to the set $G_F[t]$ defined by intersecting the sets
$\{\rho(z',z)=0\}$ over those $z\in F_{\gamma}[t]$.

\smallskip
{\it Remark.} Since there exist $\lambda: \Delta^{\kappa} \to
\Delta^n$ such that $(t,\lambda(t))\in F=F_1\cup \cdots \cup F_c$
$\forall \ t\in \Delta^{\kappa}$, there exists at least one $\gamma$
such that ${\rm gen \ rk}_{\C} (\pi|_{F_{\gamma, reg}}) =\kappa$.

\smallskip
Let $\Upsilon$ denote the dense open set of $t\in I^{\kappa}$ such
that $F_{\gamma}[t]=\emptyset$ for the $\gamma$'s with ${\rm gen \
rk}_{\C} (\pi|_{F_{\gamma,reg}}) < \kappa$. We proceed with ${\rm gen
\ rk}_{\C} (\pi|_{F_{\gamma, reg}}) = \kappa$, $\forall
\gamma=1,\ldots,c$ after forgetting other component and renumbering
the remaining ones.

Fix $F=$ a $F_{\gamma}$. Let $C=$ critical locus of
$\pi|_{F_{reg}}$. Denote $\widetilde{F}:=F_{reg}\backslash C$. It is
known that $C$ extends as a complex analytic subset of $F$ itself and
that ${\rm rk}_{\C} (\pi|_C) <\kappa$ ([C], ibidem).  Again for an
open dense set $\Upsilon$ of $t$ (still denoted by $\Upsilon$), we
have $\overline{\pi} (C) \not\ni t$ (Lemma 3.11).

Let $A_F:= \pi(\widetilde{F})$, $B_F:=\pi(\widetilde{F}) \cap
I^{\kappa}$. Clearly, $A_F$ is a nonempty subdomain of
$\Delta^{\kappa}$ (since $\widetilde{F}$ is connected).

If $B_F=\emptyset$, $\widetilde{F}\cap \pi^{-1}(t)$ makes no
contribution to $G_{F_{reg}}[t]$, if $t\in I^{\kappa}$. We can forget
those components $F$ since according to the desired conditions 1-2 of
Lemma 3.7, it is harmless to add equations to $G_{F_{reg}}[t]$ for
some other $t\in A_F$ that are close to $I^{\kappa}$ but do not belong
to $I^{\kappa}$.

Assume therefore that $B_F \neq \emptyset$.  Again, by $\pi(\Gamma
r(\lambda))=\Delta^{\kappa}$, there must exist at least one $F$ such
that $B_F \neq \emptyset$.  Let $m_1:= {\rm dim}_{\C} F$. Choose
$\underline{t} \in B_F \cap \Upsilon$, which is possible since $B_F$
is open and $\Upsilon$ is dense open, choose $\underline{\varepsilon}
>0$ with $\Delta_{\kappa}(\underline{t}, \underline{\varepsilon}) \cap
I^{\kappa} \subset \subset A_F \cap \Upsilon$.  For all $t\in
\Delta_{\kappa}(\underline{t}, \underline{\varepsilon})$,
$(\pi|_{\widetilde{F}})^{-1}(t)$ consists of finitely many
$(m_1-\kappa)$-dimensional complex submanifolds of $\widetilde{F}$,
since $\pi^{-1}(\underline{\Delta}_{\kappa}(\underline{t},
\underline{\varepsilon}))\cap C =\emptyset$, whence $\pi$ has constant
rank $\kappa$ over $\tilde{F} \cap
(\underline{\Delta}_{\kappa}(\underline{t},
\underline{\varepsilon})\times \Delta^n)$ and since \begin{equation}
(\pi|_{\widetilde{F}})^{-1}(t) \subset (\pi|_F)^{-1}(t) \end{equation}
and the latter has a finite number of connected components. This
number can only increase locally as $t$ moves. It is bounded on
$\overline{\Delta}_{\kappa}(\underline{t},
\underline{\varepsilon})\cap \overline{I}^{\kappa}$.  Hence we can
find a new $\underline{\underline{t}} \in I^{\kappa} \cap
\Delta_{\kappa}(\underline{t}, \underline{\varepsilon})$ in a
neighborhood of which this number of connected components is constant,
say in $\Delta_{\kappa}(\underline{\underline{t}},
\underline{\underline{\varepsilon}})\cap I^{\kappa}$. Denote again
simply this polydisc by $\Delta_{\kappa}(\underline{t},
\underline{\varepsilon})$.

Recall also that $\pi(C)\cap \Delta_{\kappa}(\underline{t},
\underline{\varepsilon})=\emptyset$, so $G_{F_{reg}}[t]=
G_{\widetilde{F}}[t]$, $\forall \ t\in \Delta_{\kappa}(\underline{t},
\underline{\varepsilon})$.

\smallskip\def\theequation{3.\arabic{equation}}
\setcounter{equation}{13}
\smallskip
\noindent
{\bf Lemma 3.13.} {\it Let $\underline{t}\in I^{\kappa}$ such that
there exists $\underline{\varepsilon} >0$ such that the number of
connected components of $(\pi|_{\widetilde{F}})^{-1}(t)$ is constant
equal to $\delta\in \N_*$ for all $t\in \Delta_{\kappa}(\underline{t},
\underline{\varepsilon})\cap I^{\kappa}$ and with
$\Delta_{\kappa}(\underline{t}, \underline{\varepsilon}) \cap
I^{\kappa} \subset \Upsilon$.} Then there exist holomorphic equations
$s(t,z')$ in $\Delta_{\kappa}(\underline{t},
\underline{\varepsilon})\times \Delta^{n'}$ such that

1. $G_{F_{reg}} [t]= G_{\widetilde{F}}[t] \subset \{z'\in \Delta^{n'}
\ \! {\bf :} \ \!  s(t,z')=0 \}$ \ \ $\forall \ t\in
\Delta_{\kappa}(\underline{t}, \underline{\varepsilon})$.

2. $G_{F_{reg}} [t]= G_{\widetilde{F}}[t]= \{z'\in \Delta^{n'} \ \!
{\bf :} \ \!  s(t,z')=0 \}$ \ \ $\forall \ t\in
\Delta_{\kappa}(\underline{t}, \underline{\varepsilon}) \cap
I^{\kappa}$.

\smallskip
Assume for a while that Lemma 3.13 is proved. Then Lemma 3.7 holds for
one irreducible component $F$ of the $F_{\gamma}$'s. Pick a second
component. Letting $t$ vary now in $ \Delta_{\kappa}(\underline{t},
\underline{\varepsilon})$ (instead of $\Delta^{\kappa}$), we can
repeat the above argument a finite number of steps and get Lemma 3.7
as desired. \hfill $\square$

For short, let us denote
$\underline{\Delta}_{\kappa}:=\Delta_{\kappa}(\underline{t},
\underline{\varepsilon})$.

\smallskip
{\it Proof of Lemma 3.13.} Let $D_1, \ldots,D_{\delta}$ be the
components of $(\pi|_{\widetilde{F}})^{-1}(\underline{t})$. These are
$(m_1-\kappa)$-dimensional connected complex submanifolds of
$\widetilde{F}$ (because $\pi: \widetilde{F} \to \Delta^{\kappa}$ is
submersive).  Let $p_1,\ldots,p_{\delta}\in D_1,\ldots,D_{\delta}$ be
points, let $U_1,\ldots,U_{\delta}$ be neighborhoods of
$p_1,\ldots,p_{\delta}$ in $\widetilde{F}$ with maps $\Phi_j:
\underline{\Delta}^{\kappa} \times \Delta^{m_1-\kappa} \to U_j$ such
that $\Phi_j(0\times \Delta^{m_1-\kappa})= D_j \cap U_j$,
$\Phi_j(q\times \Delta^{m_1-\kappa})$ is the fiber
$\pi^{-1}(\pi(\Phi_j(q\times 0)))\subset \tilde{F}$, $\forall \ 1\leq
j\leq \delta$, $q\in \underline{\Delta}^{\kappa}$ and such that
$\pi(\Phi_j(q\times \Delta^{m_1-\kappa}))=q$.

After all the above reductions and simplifications, we now can prove
the main step in two lemmas:

\smallskip
\noindent
{\bf Lemma 3.14.}  {\it 1.  $G_{F_{reg}}[t] \subset G_{U_1}[t] \cap
\cdots \cap G_{U_{\delta}}[t]$ \ $\forall t\in
\underline{\Delta}_{\kappa}$ and

2. $G_{F_{reg}}[t] = G_{U_1}[t] \cap \cdots \cap G_{U_{\delta}}[t]$ \
$\forall t\in \underline{\Delta}_{\kappa}\cap I^{\kappa}$.}

\smallskip\def\theequation{3.\arabic{equation}}
\setcounter{equation}{15}
\smallskip
\noindent
{\bf Lemma 3.15.} {\it Each $G_{U_j}[t]$ is equal to a set $\{ z'\in
U' : s_j(t,z')=0\}$, where $s_j=$ a finite set of holomorphic
functions.}

\smallskip
{\it Proof of Lemma 3.14.} Let $D_1[t], \ldots,D_{\delta}[t]$ denote
the connected components of $(\pi|_{\widetilde{F}})^{-1}(t)$, $t\in
\underline{\Delta}_{\kappa}\cap I^{\kappa}$.

Then $D_j[t]\cap U_j =\Phi_j(t\times \Delta^{m_1-\kappa}):=U_j[t]$ and
                      \begin{equation} G_{F_{reg}}[t] =
                      G_{\widetilde{F}}[t] = G_{{D_1}[t]}[t]\cap
                      \cdots \cap G_{D_{\delta}[t]}[t] \ \ \forall \
                      t\in \underline{\Delta}_{\kappa}\cap I^{\kappa}
\end{equation}

Now, if $\rho(z',z)=0$ $\forall \ z\in U_j[t]$, by the uniqueness
principle, then $\rho(z',z)\equiv 0$ on the connected complex manifold
$D_j[t]$, so $G_{\widetilde{F}}[t]=G_{U_1[t]}[t]\cap \cdots
G_{U_{\delta}[t]}[t]$, $\forall \ t\in \underline{\Delta}_{\kappa}\cap
I^{\kappa}$.  If $t\in \underline{\Delta}_{\kappa}\backslash
I^{\kappa}$, the cardinal of the set of connected components of
$(\pi|_{\widetilde{F}})^{-1}$ can be $>\delta$, so
$G_{\widetilde{F}}[t]$ diminishes, $G_{\widetilde{F}}[t] \subset
G_{U_1[t]}[t] \cap \cdots \cap G_{U_{\delta}[t]}[t]$.  \hfill
$\square$

\smallskip
{\it Proof of Lemma 3.15.} Let $E_j:=
\Phi_j(\underline{\Delta}_{\kappa}\times 0)$ given by a holomorphic
graph $z=\tilde{\omega} (t)$ over $\underline{\Delta}_{\kappa}$. Then
$E_j$ is a transverse manifold to the fibers of $\pi$. For each $j$,
there exist vector fields $L_1^j$, \ldots, $L_{m_1-\kappa}^j$ over
$U_j$ with holomorphic coefficients in $(t,z)$ commuting with each
other with integral manifolds $\Phi_j(t\times
\Delta^{m_1-\kappa})$. Then $\rho(z', z)=0$ $\forall \ z\in
(\pi|_{U_j})^{-1}(t)$ if and only if
$(\L^j)^{\gamma}\rho(z',z)|_{z=\tilde{\omega}(t)}=0$ $\forall \ \gamma
\in \N^{m_1-\kappa}$. Put
$s_j(t,z'):=((\L^j)^{\gamma}\rho(z',z)|_{z=\tilde{\omega}(t)})_{\gamma\in
\N^{m_1-\kappa}}$ and use noetherianity.  \hfill $\square$

\smallskip

Now, we come to the presentation of a lemma of shrinking of $\V_z'$ to
$\W_z'$ with nicer properties.

{\it Proof of Proposition 16.}

Starting with $\S_{\cal M}^1:=\{(z,\bar{w}, z') \ \! {\bf :} \ \!
(z,\bar{w}) \in {\cal M}, r(z,\bar{w},z') = 0\}$ and $\Gamma r
(f)=\{(z,\bar{w}, f(z)) \ \! {\bf :} \ \!  (z,\bar{w}) \in {\cal
M}\}\subset \S_{\cal M}^1$, we can again formalize the data as
follows. We take coordinates on ${\cal M} \cong \Delta^{\kappa}$,
$\kappa = 2m +d$.

Let $\kappa \in \N_*$, $n\in \N_*$, $J\in \N_*$, $r: \Delta^{\kappa}
\times \Delta^n \to \C^J$, $(t, z) \mapsto r(t,z)$ be a holomorphic
power series mapping converging uniformly in $(2\Delta)^{\kappa}$,
assume $r_j \in {\cal O}_{\kappa} (\Delta^{\kappa}) \times {\cal A}_n
(\Delta^n)$, let \begin{equation} S=\{(t,z)\in \Delta^{\kappa}\times
\Delta^n \ \! {\bf :} \ \!  r(t,z)=0\}.  \end{equation} Assume that
there exists $\lambda : \Delta^{\kappa} \times \Delta^n \to \Delta^n$
holomorphic, converging in $(2\Delta)^{\kappa}$ such that $\Gamma r
(\lambda) \subset S$, let $\pi: \Delta^{\kappa} \times \Delta^n \to
\Delta^n$ be the projection.

Let us inductively define a collection of $S_{\alpha}$'s, $\alpha \in
\N_*$. First $S_1= S$. Next, $S_{\alpha}= \{(t,z) \in \Delta^{\kappa}
\times \Delta^n \ \! {\bf :} \ \!  r_{\alpha} (t,z)=0\}$, $r_{\alpha}
: \Delta^{\kappa} \times \Delta^n\to \C^{J_{\alpha}}$, $J_{\alpha} \in
\N_*$, $J_{\alpha} \geq J_{\alpha-1}$, $r_{\alpha, j}= r_{\alpha- 1,
j}$ $\forall \ 1\leq j \leq J_{\alpha-1}$, $r_{\alpha, j} \in {\cal
O}_{\kappa} (\Delta^{\kappa})\times {\cal A}_n (\Delta^n)$ and $\Gamma
r(\lambda) \subset S_{\alpha}$.

The construction of $S_{\alpha+1}$ consists in forming the Jacobian
matrix of the $r_{\alpha, j}$'s with respect to $z$, $H_{\alpha}=
(\frac{\partial r_{\alpha, j}}{\partial z_k})_{1\leq k\leq n}^{1\leq j
\leq J_{\alpha}}$, in taking $(r_{\alpha+1, j})_j:=$ the collection of
all the minors $\delta_{\alpha,j}(t,z)$, $ 1\leq j\leq e_{\alpha}$, of
maximal generic rank over $\Delta^{\kappa} \times \Delta^n$ of this
matrix, where $e_{\alpha}=:J_{\alpha+1} - J_{\alpha} \in \N_*$ is the
number of such minors. Then put $(r_{\alpha+1, j})_{1\leq j \leq
J_{\alpha+1}}:= ((r_{\alpha, j})_{1\leq j \leq J_{\alpha}},
(\delta_{\alpha,j-J_{\alpha}})_{J_{\alpha}+1\leq j\leq J_{\alpha+1}}$
and put \begin{equation} S_{\alpha+1}:=\{(t,z) \in \Delta^{\kappa}
\times \Delta^n \ \! {\bf :} \ \!  r_{\alpha+1,j} (t,z)=0, 1 \leq j
\leq J_{\alpha+1}\}.  \end{equation} Of course, $r_{\alpha+1, j}\in
{\cal O}_{\kappa}(\Delta^{\kappa})\times {\cal A}_{n}(\Delta^n)$,
$\forall \ j= J_{\alpha}+1,\ldots, J_{\alpha+1}$. Also if we were
starting with $(z,\bar{w}) = t \in {\cal M}$, we would have got some
$r_{\alpha+1, j}(z,\bar{w},z')$ depending on the two variables
$(z,\bar{w})$ even if we let $(z,\bar{w})$ vary only in ${\cal M}$,
hence getting new equations like the $\tilde{r}_j$ of Proposition 16.

Then $S_{\alpha+1} \varsubsetneq S_{\alpha}$. Indeed by construction
$J_{\alpha+1} > J_{\alpha}$ and the zero-locus of equations from a
minor $\delta_{\alpha,j}$ of maximal generic rank coincides with
$S_{\alpha}$ at each point $(t_p, z_p)$ where $\delta_{\alpha}(t_p,
z_p)\neq 0$ but $S_{\alpha+1}$ does not contain $S_{\alpha}$ in a
neighborhood of such a point.

Thus there exists an integer $a\in \N_*$ such that $S_{a+1}=S_a$ and
$S_{a+1} \supset \Gamma r(\lambda)$ or $S_a \supset \Gamma r(\lambda)$
and there exists a minor $\delta_{a,j}(t,z)$ such that $\Gamma r
(\lambda) \not\subset \{\delta_{\alpha, j}=0\}$. The case $S_{a+1}=
S_a$ and $\Gamma r (\lambda) \subset S_{a+1}$ is impossible because
then ${\rm dim}_{\C} S_{a+1} \geq \kappa \geq 1$ and therefore its
minors are nontrivial which implies that $S_{a+1} \varsubsetneq S_a$
by the above remark.

Therefore $S_a\supset \Gamma r(\lambda)$ and $\Gamma r(\lambda)
\not\subset \{\delta_{a,j}=0\}$.

At each point of the Zariski open subset $\{\delta_{a,j} \neq 0\}\cap
\Gamma r(\lambda)$ of $\Gamma r (\lambda)$, locally $S_a$ is given by
equations of the form $z_2= \Phi(t,z_1)$, $(z_1, z_2)\in
\C^{n_1}\times \C^{n_1}$, $n_1+n_2 = n$, because of the constant rank
theorem. This proves $(*)$ of Proposition 16 in this context. Notice
that we make localization in a smaller open set, which is a
neighborhood of some point $(\underline{t}, \lambda(\underline{t}))\in
\Gamma r(\lambda) \cap \{\delta_{a,j} \neq 0\}$.

Next, we compute $G_F[t]$ in case $F (= S_a)$ is given by $\{(t,z)\in
\Delta^{\kappa} \times \Delta^n\ \! {\bf :} \ \!  z_2= \Phi(t,z_1)\}$
to get $(**)$. This is a particular case of Lemma 3.15: let $L=
\frac{\partial }{\partial z_1}+ \Phi_{z_1}(t,z_1) \frac{\partial
}{\partial z_2}$ be in vectorial notation the basis of vector fields
tangent to $F$. Then $\rho'(z', z_1, \Phi(t,z_1))=0$ $\forall \ z_1$
if and only if $L^{\gamma} \rho'(z', 0, \Phi (t, 0))=0$ $\forall \
\gamma \in \N^{n_1}$: these define analytic equations $s(t,z')$, which
completes the proof of Proposition 16.  Notice that we make
localization before computing $G_F[t]$: this corresponds to taking
$\W_{w, \bar{z}_1}'= \widetilde{\S}^1_{w,\bar{z}_1} \cap (U_p\times
\overline{U}_p\times U_{p'}')$ and then $r_{M'}^{U_{p'}'}(\W_{w,
\bar{z}_1}')$.  \hfill $\square$

\section{Examples}

This section is devoted to the check of examples quoted in the
introduction.  The general idea of all of these examples is to
construct $M, M', f$ with the reflection set $\S^1=\{(z,\bar{w}, z') \
\! {\bf :} \ \!  r(z,\bar{w}, z')\}$ containing two or more
irreducible components and to exploit this fact in order to exhibit
rather disharmonious phenomena about comparison between
$\X_{z,\bar{w}}'$ and $\Z_{z,\bar{w}}'$.

\def\theequation{4.\arabic{equation}} \setcounter{equation}{0}
\smallskip

{\it Check of Example 11.} Let $z\in Q_{\bar{w}}$, \begin{equation}
                      z_5= \bar{w}_5 + i z_1 \bar{w}_1.
                      \end{equation} Then $z'\in r_{M'}
                      (f(Q_{\bar{z}}))$ if and only if
                      \begin{equation} \rho'(z', (\bar{w}_1, 0, 0, 0,
                      z_5-iz_1\bar{w}_1))=0 \ \ \ \ \ \forall \
                      \bar{w}_1 \in \C, {\it i.e.}
\end{equation}
$$z_5'-[z_5-iz_1\bar{w}_1+ i\bar{w}_1^2z_3' z_4'+ iz_1' \bar{w}_1]=0 \
\ \ \ \ \forall \ \bar{w}_1 \in \C.
$$
From this follows $z_5'= z_5$, $z_1'= z_1$, $z_3' z_4'= 0$. Therefore
\begin{equation} \S^1=\{(z,\bar{w}, z')\ \! {\bf :} \ \!  z_5'=z_5,
z_1'= z_1, z_3' z_4'= 0\} \end{equation} \begin{equation}
\V_z'=\{(z_1, \zeta_2', \zeta_3', 0, z_5)\ \! {\bf :} \ \!  \zeta_2',
\zeta_3' \in \C\} \cup \{(z_1, \zeta_2', 0, \zeta_4', z_5) \ \! {\bf
:} \ \!  \zeta_2', \zeta_4' \in \C \}.  \end{equation} Now, it is
clear that $\Gamma r(f)$ is contained in \begin{equation}
\S^1_{sing}=\{(z,\bar{w}, z') \ \! {\bf :} \ \!  z_5'= z_5, z_1'= z_1,
z_3'= 0, z_4'=0\}:= \widetilde{\S}^1 \end{equation} \begin{equation}
\widetilde{\S}_{z,\bar{w}}^1=\{(z_1, \zeta_2', 0, 0, z_5) \ \! {\bf :}
\ \!  \zeta_2'\in \C\}.  \end{equation} To compute
$r_{M'}(\widetilde{\S}_{w,\bar{z}_1}^1)$, we write \begin{equation}
\rho'(w', \bar{w}_1, \bar{\zeta}_2', 0, 0, \bar{w}_5)= 0 \ \ \ \forall
\ \bar{\zeta}_2', {\it i.e.}
\end{equation}
$$w_5'- \bar{w}_5 - i[ w_1' \bar{w}_1 + \bar{w}_1^2 w_3' w_4'+
\bar{\zeta'}_2^2 {w'}_3^2 ] =0 \ \ \ \forall \ \bar{\zeta}_2'.
           $$
From this follows $w_3'= 0$, $w_5'= \bar{w}_5 +iw_1'
                      \bar{w}_1$. Therefore \begin{equation} r_{M'}
                      (\widetilde{\S}_{w,\bar{z}_1}^1)= \{(w_1', w_2',
                      0, w_4', \bar{w}_5+iw_1'\bar{w}_1)\ \! {\bf :} \
                      \!  w_1', w_2', w_4' \in \C\} \end{equation} and
                      \begin{equation} \widetilde{\S}_{z,\bar{w}}^1
                      \cap \widetilde{\S}_{w, \bar{z}_1}^1= \{(z_1,
                      \zeta_2', 0,0, z_5) \ \! {\bf :} \ \!  \zeta_2'
                      \in \C\}.  \end{equation} In conclusion, ${\rm
                      dim}_{f(z)} \Z_{z,\bar{w}}'= $ $\forall \ z$. On
                      the other hand, \begin{equation} \V_w'=\{(w_1,
                      \zeta_2', \zeta_3', 0, w_5) \ \! {\bf :} \ \!
                      \zeta_2', \zeta_3' \in \C\} \cup \{(w_1,
                      \zeta_2', 0, \zeta_4', w_5) \ \! {\bf :} \ \!
                      \zeta_2', \zeta_3'\in \C \} \end{equation} and
                      the equations of $r_{M'}(\V_w')$ are given by
                      \begin{equation} \rho'(w', (\bar{w}_1,
                      \bar{\zeta}_2', \bar{\zeta}_3', 0,
                      \bar{w}_5'))=0 \ \ \ \ \ \forall \
                      \bar{\zeta'}_2 \forall \ \bar{\zeta}_3', {\it
                      i.e.}
\end{equation}
$$w_5' -[\bar{w}_5+ i[w_1' \bar{w}_1+ \bar{w}_1^2 w_3'w_4'+ w_3'
\bar{\zeta'}_2^2 +{w'}_2^2\bar{\zeta'}^2_3+
{w'}_4^3\bar{\zeta'}_3^3]]=0 \ \ \ \ \ \forall \ \bar{\zeta'}_2
\forall \ \bar{\zeta}_3'
$$
                      \begin{equation} \rho'(w', (\bar{w}_1,
                      \bar{\zeta}_2', 0, \bar{\zeta}_4', \bar{w}_5))=0
                      \ \ \ \ \ \forall \bar{\zeta'}_2 \forall
                      \bar{\zeta'}_4, {\it i.e.}
\end{equation}
$$w_5'- [\bar{w}_5 + i[w_1' \bar{w}_1+ \bar{w}_1^2 w_3' w_4'+ {w'}^2_3
\bar{\zeta'}_2^2 + \bar{\zeta'}_4^3 {w'}_3^3]]=0 \ \ \ \ \ \forall
\bar{\zeta'}_2 \forall \bar{\zeta'}_4.
$$
From (4.11) we deduce $w_3'=0$, $w_2'=0$, $w_4'=0$, $w_5'= \bar{w}_5+
i w_1' \bar{w}_1$. From (4.12) we deduce $w_3'= 0$, $w_5'= \bar{w}_5+
iw_1' \bar{w}_1$. Therefore \begin{equation} r_{M'} (\V_w')= \{ (w_1',
0,0,0, \bar{w}_5+ iw_1' \bar{w}_1) \ \! {\bf :} \ \!  w_1' \in \C\} =
f(Q_{\bar{w}}) \end{equation} and finally \begin{equation} \V_z' \cap
r_{M'} (\V_w')= \{(z_1,0,0,0, z_5)\}= \{f(z)\}.  \end{equation} In
conclusion, ${\rm dim}_{f(z)} \X_{z,\bar{w}}'=0$. This completes
Example 11. \hfill $\square$

\smallskip
{\it Check of Example 12.} Here, if $z_4= \bar{w}_4+i\bar{w}_1 z_1$,
\begin{equation} \S^1= \{(z, \bar{w}, z') \ \! {\bf :} \ \!  z_4'=
z_4, z_1'= z_1, z_2' z_3'= 0\}.  \end{equation} Now, it is clear that
$\Gamma r(f)$ is contained in ${\S}^1_{sing}$ \begin{equation}
\S^1_{sing}= \{(z,\bar{w}, z') \ \! {\bf :} \ \!  z_4'=z_4, z_1'= z_1,
z_2'= 0, z_3'= 0\}.  \end{equation} Whence \begin{equation}
\widetilde{\S}^1_{z,\bar{w}}=\{(z_1,0,0,0, z_5)\}= \W_{z,\bar{w}}'
\end{equation} and finally \begin{equation} {\rm dim}_{f(z)}
\W_{z,\bar{w}}' = {\rm dim}_{f(z)} \Z_{z,\bar{w}}'= 0 \ \ \ \ \
\forall \ z.  \end{equation} On the other hand, \begin{equation}
\S_{z,\bar{w}}^1= \{(z_1, \zeta_2', 0, z_4) \ \! {\bf :} \ \!
\zeta_2'\in \C\} \cup \{(z_1, 0, \zeta_3', z_4) \ \! {\bf :} \ \!
\zeta_2' \in \C\} = \V_z' \end{equation} and the equations of
$r_{M'}(\V_w')$ are given by \begin{equation} \rho'(w', (\bar{w}_1, 0,
\bar{\zeta'}_3, \bar{w}_4))= w_4' - [\bar{w}_4+ i [w_1' \bar{w}_1+
\bar{w}_1^2 w_2' w_3']] = 0 \ \ \ \ \ \forall \ \bar{\zeta}_3'
\end{equation} \begin{equation} \rho'(w', (\bar{w}_1, \bar{\zeta'}_2,
0, \bar{w}_4))= w_4'- [\bar{w}_4 + i [w_1' \bar{w}_1 + \bar{w}_1^2
w_2' w_3' ]]= 0 \ \ \ \ \ \forall \bar{\zeta}_2'.  \end{equation} It
follows only the equation $w_4'= \bar{w}_4 + i[w_1'\bar{w}_1 +
\bar{w}_1^2 w_2' w_3']$. Therefore \begin{equation} r_{M'} (\V_w')=
\{(w_1', w_2', w_3', \bar{w}_4 + i[w_1' \bar{w}_1+ \bar{w}_1^2 w_2'
w_3']) \ \! {\bf :} \ \!  w_1', w_2', w_3' \in \C\} \end{equation}
\begin{equation} \V_z' \cap r_{M'} (\V_w') = \{(z_1, w_2', 0, z_4) \
\! {\bf :} \ \!  w_2' \in \C\} \cup \{(z_1, 0, w_3', z_4) \ \! {\bf :}
\ \!  w_3' \in \C\}= \V_z', \end{equation} whence ${\rm dim}_{f(z)}
\X_{z,\bar{w}}' = 1$ $\forall \ z$. Example 12 is complete. \hfill
$\square$

\smallskip
{\it Check of Example 13.} Let us establish:

\smallskip
{\it The function ${\cal M} \ni (z,\bar{w}) \mapsto {\rm dim}_{f(z)}
\X_{z,\bar{w}}'\in \N$ is neither upper semi continuous nor lower semi
continuous in general.}

\smallskip
{\it Proof.} First, whenever ${\rm dim}_{f(0)} r_{M'} (f(Q_0))=0$,
then ${\rm dim}_{f(z)} r_{M'} (f(Q_{\bar{z}})) = 0$ too for $z\in
{\cal V}_{\C^n} (0)$ because of Lemma 2.3 and so there exists $V=
{\cal V}_{\C^n}(0)$ such that ${\rm dim}_{f(z)} \X_{z,\bar{w}}'=0$
$\forall \ z, w \in U$, $z\in Q_{\bar{w}}$. For instance, $M=M'$, $f=
{\rm Id}$, $M=\{z_2= \bar{z}_2+iz_1\bar{z}_1\}$.

Therefore $(z,\bar{w})\mapsto {\rm dim}_{f(z)} \X_{z,\bar{w}}'$ could
be continuous.

{\it This is false.} Indeed, let $M= M'= \{z_3= \bar{z}_3 + i
z_1\bar{z}_1 (1+ z_2 \bar{z}_2) \}\subset \C^3$, $f= {\rm Id}$. First,
$M$ is Levi-nondegenerate at every point of $\C^3 \backslash \{z_1=
0\}$, so ${\rm dim}_{f(p)} \X_{p, \bar{p}}'= 0$ at those points. Let
$Q_0= \{(z_1, z_2, 0) \ \! {\bf :} \ \!  z_1, z_2 \in \C\}$,
$r_{M'}(Q_0)=\{(0, z_2, 0) \ \! {\bf :} \ \!  z_2 \in \C\}= \{q\in
\C^3 \ \! {\bf :} \ \! Q_{\bar{q}}= Q_0\}$, so $r_{M'} (Q_0)= Q_0=
\{(z_1, z_2, 0) \ \! {\bf :} \ \!  z_1, z_2 \in \C\}$, so $r_{M'}
(Q_0) \cap r_{M'}^2(Q_0)= \{(0, z_2, 0)\ \! {\bf :} \ \!  z_2 \in
\C\}$ has dimension 1.

Therefore $(z,\bar{w})\mapsto {\rm dim}_{f(z)} \X_{z,\bar{w}}'$ can be
at best upper semi-continuous.

{\it This is false.} Indeed, let $M= \{(z_1, z_4) \in \C^2 \ \! {\bf
:} \ \!  z_4 =\bar{z}_4 + i z_1 \bar{z}_1\}$, let \begin{equation} M'
= z_4'= \bar{z}_4'+ i z_1'\bar{z}_1'+ i z_3' \bar{z}_2'+ i \bar{z}_3'
z_2'+ i {z'}_1^2 \bar{z}_3' \bar{z}_2'+ i\bar{z'}_1^2 z_3' z_2'
\end{equation} and \begin{equation} f(z_1, z_4)= (z_1, z_4 \sin^3 z_1,
0, z_4), \ \ \ \ \ f(M) \subset M'.  \end{equation} Then $Q_0 =
\{(z_1, 0)\}$, $f(Q_0) = \{(z_1, 0,0,0)\}$. We claim that

1. $r_{M'}(f(Q_0)) \cap r_{M'} (f(Q_0)) = \{(0,0,0,0)\}$ and

2. $r_{M'}(f(Q_{\bar{z}})) \cap r_{M'}^2(f(Q_{\bar{w}}))=\{ (z_1,
z_2', 0, z_4) \ \! {\bf :} \ \!  z_2'\in \C\}$ $\forall \ z, w$, $z\in
Q_{\bar{w}}$, $z\neq 0$.

\noindent
This will show that $(z,\bar{w})\mapsto {\rm dim}_{f(z)}
\X_{z,\bar{w}}'$ cannot be upper semicontinuous in general.

Indeed,
\begin{equation}
r_{M'}(f(Q_0))= \{(0, z_2, 0,0) \} \cup \{(0,0,z_3,0)\} \ \ {\rm and}
\ \ r_{M'}^2(f(Q_0))= \{(z_1,0,0,0)\},
\end{equation} 
so 1 holds.

Let $(z_1,z_4)\neq (0,0)$, let $z\in Q_{\bar{w}}$, $w=(w_1,w_4)$,
$z_4=\bar{w}_4+iz_1\bar{w}_1$. Then $f(Q_{\bar{z}})=\{(w_1,
(\bar{z}_4+iw_1\bar{z}_1) \sin^3 w_1, 0, \bar{z}_4+iw_1\bar{z}_1)\ \!
{\bf :} \ \!  w_1\in \C\}$. By definition, $r_{M'}
(f(Q_{\bar{z}}))=\{z'\ \! {\bf :} \ \!  \rho'(z', \bar{f}(\bar{w}))=0
\ \forall \ w\in Q_{\bar{z}}\}$. Write \begin{equation} \rho'(z',
\bar{f}(\bar{w}))=0 \ \ \ \ \ \forall \ \bar{w}_1, \ {\it i.e.}
\end{equation}
$$z_4'- [z_4- i \bar{w}_1 z_1+ iz_1' \bar{w}_1+ iz_3'(z_4- i\bar{w}_1
z_1) \sin^3 \bar{w}_1+ iz_3' z_2' \bar{w}_1^2]= 0 \ \ \ \ \ \forall \
\bar{w}_1.
$$
We deduce equations $z_4'= z_4$, $z_1'= z_1$, $z_3'z_2'= 0$, $z_3'
z_4=0$, $z_3' z_1= 0$. Therefore if $(z_1, z_4)\neq (0,0)$, then
\begin{equation} r_{M'}(f(Q_{\bar{z}}))= \{(z_1, z_2', 0, z_4) \ \!
{\bf :} \ \!  z_2' \in \C\}.  \end{equation} Next,
$r_{M'}^2(f(Q_{\bar{w}}))$ is given by \begin{equation} \rho'(z',
(\bar{w}_1, \bar{\zeta}_2', 0, \bar{w}_4))= z_4' - [\bar{w}_4 + i z_1'
\bar{w}_1 + i z_3' z_2' \bar{w}_1^2 + i z_3' \bar{\zeta}_2']= 0 \ \ \
\ \ \forall \ \bar{\zeta}_2'.  \end{equation} We deduce equations
$z_3'= 0$, $z_4'= \bar{w}_4 + i z_1' \bar{w}_1$, so \begin{equation}
r_{M'}^2(f(Q_{\bar{w}}))= \{ (z_1', z_2', 0, \bar{w}_4+ i z_1'
\bar{w}_1) \ \! {\bf :} \ \!  z_1', z_2' \in \C\}, \ \ \ \ \ (z_1,
z_4) \neq (0,0).  \end{equation} Finally for such $(z_1, z_4) \neq
(0,0)$, \begin{equation} r_{M'} (f(Q_{\bar{z}})) \cap r_{M'}^2
(f(Q_{\bar{w}})) = \{ (z_1, z_2', 0, z_4) \ \! {\bf :} \ \!  z_2' \in
\C\}, \end{equation} which shows that 2 above holds. This completes
Example 13. \hfill $\square$

\smallskip
Example 13 already shows that $\X_{z,\bar{w}}'$ is not analytically
parametrized by $(z,\bar{w})$.  Example 14 also provides a
supplementary reason.

\smallskip
{\it Check of Example 14.} First, let us take in (3.4): $\kappa=2$,
$n=3$, $\lambda(t_1, t_2)= (t_1, t_2, 0)$, \begin{equation} F=\{(t, z)
\in \Delta^{2} \times \Delta^3 \ \! {\bf :} \ \!  (t_1z_3- t_2^2)
(z_1- t_1) =0, (t_1 z_3 - t_2^2) (z_2- t_2)=0,
\end{equation}
$$(t_1 z_3- t_2^2) z_3 = 0\} = F_1 \cup F_2 = \Gamma r (\lambda) \cup
\{ t_1 z_3- t_2^2= 0\}.
$$
Then the fibers $F_2[t] =\emptyset$ if $|t_2^2|\geq |t_1|$, say if
$t\in T_c := \Delta^2 \cap \{|t_2^2| \geq |t_1|\}$, and $F_2[t]=
\{(z_1, z_2, z_3) \in \Delta^3\ \! {\bf :} \ \!  z_3= t_2^2 / t_1\}$
if $t\in T_0:= \Delta^2 \backslash T_c$. Clearly \begin{equation}
G_F[t] = \{z'\in \Delta^{n'} \ \! {\bf :} \ \!  \rho(z', (t_1, t_2,
0))= 0\} \ \ \ \ \ \forall \ t\in \Delta^2 \end{equation} and
\begin{equation} G_{F_2}[t]= \{z'\in \Delta^{n'}\} \end{equation} if
$t\in T_c$ and \begin{equation} G_{F_2}[t]= \{z'\in \Delta^{n'} \ \!
{\bf :} \ \!  (\partial_{z_1}^{k_1} \partial_{z_2}^{k_2} \rho)(z',
0,0, t_2^2 / t_1)= 0 \ \forall \ k_1 \forall \ k_2 \in \N\}
\end{equation} if $t\in T_0$. The border equals $\overline{T}_0 \cap
T_c=\{t\in \Delta^2 \ \! {\bf :} \ \!  |t_1| = |t_2^2|\}$. It is real
analytic, not complex.

Next, we build a mapping on the basis of this example. Let $n=3$,
$n'=4$, $M: z_4=\bar{z}_4+iz_1\bar{z}_1+ iz_2\bar{z}_2$,
\begin{equation} M': \ \ z_4'= \bar{z}_4'+ i[{z'}_1^2
\bar{z}_3'(\bar{z'}_2^2 - \bar{z}_1' \bar{z}_3') + \bar{z'}_1^2
z_3'({z'}_2^2 - z_1' z_3')+ z_1'\bar{z}_1'+ z_2'\bar{z}_2'],
\end{equation} $f(z_1, z_2, z_4)= (z_1, z_2, 0, z_4)$. Then one can
check that the equations of $\S_{z,\bar{w}}'= \V_z'$ are: $z_4'= z_4$,
$z_1'= z_1$, $z_2'=z_2$, $z_3'(z_2^2- z_1z_3)= 0$, from which Example
14 follows. \hfill $\square$

\smallskip
{\it Check of Example 15.}  Consider $f: \C^2 \ni (z_1, z_4) \mapsto
(z_1, 0,0,z_4) \in \C^4$, $M: z_4= \bar{z}_4 + iz_1\bar{z}_1$ and
\begin{equation} M': \ z_4'= \bar{z}_4' + iz_1' \bar{z}_1' + i
{z'}_1^2 (1+ \bar{z}_3')\bar{z}_3' + i \bar{z'}_1^2(1+ z_3')z_3'+ i
z_2' z_3' \bar{z'}_2^2 + i \bar{z}_2' \bar{z}_3' {z'}_2^2.
\end{equation} Identify $M$ with $f(M)= \{(z_1, 0,0,z_4) \ \! {\bf :}
\ \!  (z_1, z_4) \in M\}$. Take $U= \Delta^2 \cong \Delta \times 0
\times 0 \times \Delta$, $U'= \Delta^4$. We will first check 1 and 2
for $z=\bar{w}= 0\in \C^2$.

First, compute $r_{M'}(Q_0)$ by writing $\rho'(z', \bar{f}(\bar{z}_1,
0))= z_4'- [ i\bar{z}_1 z_1' + i \bar{z}_1^2 (1+ z_3') z_3']$ so that
equations of $r_{M'} (Q_0)$ are $z_4'= 0$, $z_1'= 0$, $(1+z_3')z_3'=
0$, whence \begin{equation} r_{M'}^{U'}(Q_0)= \{(0, z_2',0,0) \ \!
{\bf :} \ \!  z_2' \in \Delta \} \end{equation} \begin{equation}
r_{M'}^{\C^{n'}} (Q_0)= \{(0, z_2', 0, 0) \ \! {\bf :} \ \!  z_2' \in
\C\} \cup \{(0, z_2', 1, 0) \ \! {\bf :} \ \!  z_2' \in \C\}:= \A_0^1
\cup \A_0^2.  \end{equation} To compute $(r_{M'}^{U'})^2(Q_0)$, write
$\rho'(z', (0, \bar{w}_2', 0,0))= z_4'- i[\bar{w'}_2^2 z_2' z_3']$ so
that equations of $(r_{M'}^{U'})^2(Q_0)$ are $z_4'= 0$, $z_2' z_3'=
0$, whence \begin{equation} (r_{M'}^{U'})^2(Q_0)=\{(z_1', z_2', 0,0) \
\! {\bf :} \ \!  z_1', z_2' \in \Delta\} \cup \{(z_1', 0, z_3', 0) \
\! {\bf :} \ \!  z_1', z_3'\in \Delta\} \end{equation}
\begin{equation} (r_{M'}^{U'})(Q_0) \cap (r_{M'}^{U'})^2(Q_0)=\{ (0,
z_2', 0,0) \ \! {\bf :} \ \!  z_2' \in \Delta\}.  \end{equation} On
the other hand, $(r_{M'}^{\C^{n'}})^2(Q_0)= r_{M'}^{\C^{n'}}(\A_0^1)
\cap r_{M'}^{\C^{n'}}(\A_0^2)$, where as above \begin{equation}
r_{M'}^{\C^{n'}}(\A_0^1)=\{(z_1', z_2', 0,0) \ \! {\bf :} \ \!  z_1',
z_2' \in \C\} \cup \{(z_1', 0, z_3', 0) \ \! {\bf :} \ \!  z_1',
z_3'\in \C\}.  \end{equation} To compute $r_{M'}^{\C^{n'}}(\A_0^2)$,
write $\rho'(z', (0, \bar{w}_2', 1, 0))= z_4' - i[\bar{w'}_2^2
z_2'z_3' + \bar{w}_2' {z'}_2^2]= 0 $ $\forall \ \bar{w}_2'$ so that
its equations are $z_4'=0$, ${z'}_2^2=0$, $z_2' z_3'=0$, whence
$r_{M'}^{\C^{n'}}(\A_0^2)=\{(z_1', 0, z_3', 0) \ \! {\bf :} \ \!  z_1'
\in \C, z_2' \in \C\}$ and \begin{equation} r_{M'}^{\C^{n'}}(Q_0) \cap
(r_{M'}^{\C^{n'}})^2(Q_0)= \{(0,0,0,0)\} \cup \{(0,0,1,0)\}.
\end{equation} In conclusion, for $z=\bar{w}=0\in M$,
$r_{M'}^{\C^{n'}}(Q_0) \cap (r_{M'}^{\C^{n'}})^2(Q_0)$ is finite
whereas ${\rm dim}_0 [ r_{M'}^{U'}(Q_0) \cap (r_{M'}^{U'})^2(Q_0)]=
1$. Now, let $z\in Q_{\bar{w}}$, $z,w\in \Delta \times 0 \times 0
\times \Delta$, $(z_1,0,0,z_4)\in
r_{M'}^{U'}(Q_{\bar{z}_1,0,0,\bar{z}_4})\cap
(r_{M'}^{U'})^2(Q_{\bar{w}_1,0,0,\bar{w}_4})$. As above,
\begin{equation} r_{M'}^{U'} (Q_{\bar{z}})=\{(z_1, z_2', 0, z_4) \ \!
{\bf :} \ \!  z_2'\in \Delta\} \end{equation} \begin{equation}
r_{M'}^{\C^{n'}} (Q_{\bar{z}})= \{(z_1,z_2', 0, z_4) \ \! {\bf :} \ \!
z_2' \in \C\} \cup \{(z_1, z_2', 1, z_4) \ \! {\bf :} \ \!  z_2' \in
\C\} := \A_z^1 \cup \A_z^2.  \end{equation} To compute
$(r_{M'}^{U'})^2 (Q_{\bar{w}})= \{w'\in \Delta^4 \ \! {\bf :} \ \!
\rho'(w', \bar{z}') = 0 \ \forall \ z'\in r_{M'}^{U'}(Q_{\bar{w}}\})$,
write first \begin{equation} \rho'(w', (\bar{w}_1, \bar{z}_2', 0,
\bar{w}_4))= w_4'-[w_4+ i w_1' \bar{w}_1+i \bar{w}_1^2 (1+ w_3') w_3'
+ i w_2' w_3' \bar{z'}_2^2] \end{equation} whence \begin{equation}
(r_{M'}^{U'})^2(Q_{\bar{w}})=\{(w_1', w_2', 0, \bar{w}_4+i
w_1'\bar{w}_1) \ \! {\bf :} \ \!  w_1', w_2'\in \Delta, |\bar{w}_4 + i
w_1' \bar{w}_1 | < 1 \} \cup \end{equation}
$$ \ \ \ \ \ \ \ \ \ \cup \{(w_1', 0, w_3', \bar{w}_4+ i w_1'
\bar{w}_1) \ \! {\bf :} \ \!  w_1', w_2' \in \Delta, |\bar{w}_4 + i
w_1' \bar{w}_1| < 1 \}.$$ \begin{equation} r_{M'}^{U'}(Q_{\bar{z}})
\cap (r_{M'}^{U'})^2(Q_{\bar{w}})= \{(z_1, w_2', 0, z_4) \ \! {\bf :}
\ \! w_2' \in \Delta\}.  \end{equation} On the other hand analogously
\begin{equation} r_{M'}^{\C^{n'}}(\A_w^1)=\{(z_1, w_2', 0, z_4) \ \!
{\bf :} \ \!  w_2' \in \C\}.  \end{equation} To compute
$r_{M'}^{\C^{n'}}(\A_w^2)= \{w'\in \C^4 \ \! {\bf :} \ \! \rho'(w',
\bar{z}')= 0 \ \forall \ z'\in \A_w^2\}$, write \begin{equation}
\rho'(w', (w_1, z_2', 1, w_4))=
\end{equation}
$$=w_4'-[w_4+i w_1' \bar{w}_1 + i\bar{w}_1^2 (1+ w_3') w_3' + i
w_2'w_3' \bar{z'}_2^2+ i {w'}_2^2 \bar{z}_2']
$$
whence \begin{equation} r_{M'}^{\C^{n'}} (\A_w^2)=\{(w_1',0,w_3',
                      \bar{w}_4+iw_1'\bar{w}_1) \ \! {\bf :} \ \!
                      w_1', w_2' \in \C\} \end{equation}
                      \begin{equation} r_{M'}^{\C^{n'}}(Q_{\bar{z}})
                      \cap (r_{M'}^{\C^{n'}})^2 (\A_w^1) \cap
                      (r_{M'}^{\C^{n'}})^2 (\A_w^2)=\{(z_1,0,0,z_4)\}
                      \cup \{(z_1,0,1,z_4\}.  \end{equation} This
                      completes Example 15. \hfill $\square$

\smallskip
\smallskip
\noindent
{\bf Closing remark.} During the completion of his work about
algebraicity of holomorphic mappings between real algebraic CR
manifolds [ME99a], the author had to study the preprint ($3^{rd}$
version 29pp; $2^{nd}$ version 18pp: spring 1998; $1^{st}$: e-print:
{\sf http://xxx.lanl.gov/abs/math.CV/9801040}; 1998 11pp)
``Algebraicity of local holomorphisms between real algebraic
submanifolds of complex spaces'' by D. Zaitsev, submitted for
publication near august 1998 in which Theorem 1' here is claimed.
While studying this paper, Damour (see [D99]) observed an incorrection
in the Lemma 4.2 in [Z98].  It appears that what is actually proved by
Zaitsev is that
$$ {\rm dim}_{f(p)} \Z_{p,\bar{p}}'=0 \ \ \forall \ p\in M\cap {\cal
V}_{\C^n}(0) \ \Rightarrow \ f \ {\rm is \ algebraic}.$$ Our examples
exhibit several phenomena which show that this statement is not
equivalent to Theorem 1': there is a serious difference between
$\Z_{z,\bar{w}}'$ and $\X_{z,\bar{w}}'$. They also show that it is
more natural to proceed as in [ME99a] by choosing at the first steps a
{\it minimal for inclusion} real algebraic CR manifold $M''$ which
satisfies $f(M)\subset M'' \subset M'$. This crucial step explains all
our examples after careful examination.  Finally, in [ME99a] the
conditions given for algebraicity are necessary and sufficient whereas
in [Z98] they are sufficient only.

\end{document}